\documentclass[11pt]{article}
\usepackage{epsfig}
\usepackage{graphicx}
\usepackage{amssymb,bm}
\usepackage{mathrsfs}
\usepackage{amsthm}

\font\absfont cmss10 at 8pt
\renewenvironment{abstract}
  {\list{}{\listparindent 0.2cm%
  \baselineskip 0pt
  \setlength{\leftmargin}{1.5cm}
  \setlength{\rightmargin}{1cm}
  }%
  \item\relax \hskip -10pt {\sc Abstract.}\ \absfont}
  {\endlist}

\DeclareMathSymbol{\varGamma}{\mathord}{letters}{"00}
\DeclareMathSymbol{\varDelta}{\mathord}{letters}{"01}
\DeclareMathSymbol{\varTheta}{\mathord}{letters}{"02}
\DeclareMathSymbol{\varLambda}{\mathord}{letters}{"03}
\DeclareMathSymbol{\varXi}{\mathord}{letters}{"04}
\DeclareMathSymbol{\varPi}{\mathord}{letters}{"05}
\DeclareMathSymbol{\varSigma}{\mathord}{letters}{"06}
\DeclareMathSymbol{\varUpsilon}{\mathord}{letters}{"07}
\DeclareMathSymbol{\varPhi}{\mathord}{letters}{"08}
\DeclareMathSymbol{\varPsi}{\mathord}{letters}{"09}
\DeclareMathSymbol{\varOmega}{\mathord}{letters}{"0A}

\theoremstyle{plain}
\newtheorem{theorem}{Theorem}

\newtheorem{lemma}{Lemma}[section]
\newtheorem{corollary}[lemma]{Corollary}
\newtheorem{proposition}[lemma]{Proposition}

\newtheorem{definition}[lemma]{Definition}
\newtheorem{assumption}{Assumption}
\theoremstyle{remark}
\newtheorem{remark}[lemma]{Remark}

\makeatletter\@addtoreset{equation}{section}
\makeatother

\newcommand{\C}{\mathbb{C}}

\newcommand{\R}{\mathbb{R}}
\newcommand{\Z}{\mathbb{Z}}
\newcommand{\p}{\partial }

\newcommand{\Abs}[1]{\left\vert#1\right\vert}
\newcommand{\abs}[1]{\vert #1\vert}

\newcommand{\norm}[1]{\Vert #1 \Vert}
\newcommand{\at}[1]{\vert\sb{#1}}
\newcommand{\At}[1]{\Big\vert\sb{#1}}
\newcommand\const{\mathop{\rm const}}

\newcommand\rank{\mathop{\rm rank}}
\newcommand\range{\mathop{\rm Range}}
\newcommand\image{\mathop{\rm Image}}

\newcommand\sothat{{\rm :}\ }

\newcommand{\fra}[2]{#1/#2}

\newfam\bifam
\font\tenbi=cmmib10 scaled \magstep1
\font\sevenbi=cmmib10 at 11pt
\font\fivebi=cmmib10 at 6pt
\textfont\bifam=\tenbi
\scriptfont\bifam=\sevenbi
\scriptscriptfont\bifam=\fivebi

\DeclareSymbolFont{EUF}{U}{euf}{m}{n}
\SetSymbolFont{EUF}{bold}{U}{euf}{b}{n}
\DeclareSymbolFontAlphabet{\euf}{EUF}

\newcommand{\os}{\mathop{\rm o}}
\newcommand{\ol}{\mathop{\rm O}}

\begin{document}

\title{\bf Nonlinear instability of a critical traveling wave
in the generalized Korteweg -- de Vries equation}

\author{
{\sc Andrew Comech}\footnote{A.C. was partially supported by
Max-Planck Institute for Mathematics in the Sciences (Leipzig)
and by the NSF Grants DMS-0434698 and DMS-0621257.
}
\\
{\small\it
Department of Mathematics, Texas A\&M University, College Station, TX 77843, USA}
\medskip
\\
{\sc Scipio Cuccagna\footnote{S.C. was fully supported by a special
grant of the Italian Ministry of Education, University and Research.}}
\\
{\small\it DISMI, University of Modena and Reggio Emilia,
Reggio Emilia 42100 Italy}
\medskip
\\
{\sc Dmitry E. Pelinovsky}
\\
{\small\it Department of Mathematics,
McMaster University, Hamilton, ON L8S\,4K1, Canada}
}

\date{\today}
\maketitle

\begin{abstract}
We prove the instability of a ``critical'' solitary wave
of the generalized Korteweg -- de Vries equation,
the one
with the speed
at the border between the stability and instability regions.
The instability mechanism involved is ``purely nonlinear'',
in the sense that the linearization at a critical soliton
does not have eigenvalues with positive real part.
We prove that critical solitons correspond generally
to the saddle-node bifurcation of two branches of solitons.
\end{abstract}

\section{Introduction and main results}\label{sect-introduction}

We consider the generalized Korteweg -- de Vries equation in
one dimension,
\begin{equation}
\label{gkdv}
\p\sb{t}{\bm u}
=\p\sb{x}\left(-\p\sb{x}^2 {\bm u}+f({\bm u})\right),
\qquad
{\bm u}={\bm u}(x,t)\in\R,\quad x\in\R,
\end{equation}
where
$f\in C\sp\infty(\R)$
is a real-valued function that satisfies
\begin{equation}\label{f-p-zero}
f(0)=f'(0)=0.
\end{equation}
Depending on the nonlinearity $f$, equation (\ref{gkdv}) may
admit solitary wave solutions, or solitons,
of the form
${\bm u}(x,t)=\bm\phi\sb{c}(x-ct)$.
Generically, solitons exist for speeds $c$
from (finite or infinite) intervals of a real line.
For a particular nonlinearity $f$,
solitons with certain speeds are (orbitally) stable with
respect to the perturbations of the initial data,
while others are linearly
(and also dynamically)
unstable.
We will study the stability of the critical solitons,
the ones with the speeds $c$ on the border of stability
and instability regions.
These solitons are no longer linearly unstable.
Still, we will prove their instability,
which is the consequence of the higher
algebraic multiplicity of the zero eigenvalue
of the linearized system.

When
$f({\bm u})=-3{\bm u}^2$,
(\ref{gkdv})
turns into the classical Korteweg -- de Vries (KdV) equation
\begin{equation}\label{kdv}
\p\sb{t}{\bm u}+\p\sb{x}\sp{3}{\bm u}+6{\bm u}\p\sb{x}{\bm u}=0
\end{equation}
which is well-known to have solitary-wave solutions, or solitons,
\[
{\bm u}\sb{c}(x,t)
=\bm\phi\sb{c}(x-ct)
=\frac{c}{2\cosh\sp{2}\left(\frac{\sqrt{c}}{2}(x-ct)\right)},
\qquad c>0.
\]
For $f({\bm u})=-{\bm u}^p$, $p>1$,
we obtain the family of
generalized KdV equations
(also known as gKdV-$k$ with $k=p-1$)
that have the form
\begin{equation}\label{kdvk}
\p\sb{t}{\bm u}+\p\sb{x}\sp{3}{\bm u}+\p\sb{x}({\bm u}^{p})=0.
\end{equation}
They also have solitary wave solutions.
All solitary waves of the
classical KdV equation and of the subcritical generalized KdV
equations ($1<p<5$) are orbitally stable; see \cite{MR0338584},
\cite{MR0386438}, \cite{MR886343}, \cite{MR887857}.
Orbital stability is defined in the following sense:

\begin{definition}\label{def-stability}
The traveling wave $\bm\phi\sb{c}(x-ct)$ is said
to be orbitally stable if for any $\epsilon>0$
there exists $\delta>0$ so that
for any $\bm u\sb 0$
with $\norm{\bm u\sb 0-\bm\phi\sb{c}}\sb{H\sp 1}\le\delta$
there is a solution ${\bm u}(t)$ with ${\bm u}(0)={\bm u}\sb 0$,
defined for all $t\ge 0$, such that
\[
\sup\sb{t\ge 0}\inf\sb{s\in\R}
\norm{{\bm u}(x,t)-\bm\phi(x-s)}\sb{H\sp 1}<\epsilon,
\]
where $H\sp 1=H\sp 1(\R)$ is the standard Sobolev space.
Otherwise the traveling wave is said to be unstable.
\end{definition}

Equation (\ref{gkdv})
is a Hamiltonian system,
with the Hamiltonian functional
\begin{equation}
\label{hamiltonian}
E({\bm u})
=\int\limits\sb{\R}\left(
\frac{1}{2}(\p\sb{x}{\bm u})^2+F({\bm u})
\right)\,dx,
\end{equation}
with $F({\bm u})$  the antiderivative of $f({\bm u})$
such that $F(0)=0$.
There are two more invariants of motion:
the mass
\begin{equation}\label{mass}
I({\bm u})=\int\limits\sb{\R} {\bm u}\,dx
\end{equation}
and the momentum
\begin{equation}\label{momentum}
\mathscr{N}({\bm u})
=\int\limits\sb{\R}\frac{1}{2}{\bm u}^2\,dx.
\end{equation}

\begin{assumption}\label{assumption-soliton}
There is an open set $\varSigma\subset\R\sb{+}$
so that for $c\in\varSigma$
the equation
$-c\bm\phi\sb{c}=-\bm\phi\sb{c}''+f(\bm\phi\sb{c})$
has a unique solution $\bm\phi\sb{c}(x)\in H\sp\infty(\R)$
such that
$\bm\phi\sb{c}(x)>0$,
$\bm\phi\sb{c}(-x)=\bm\phi\sb{c}(x)$,
$\lim\sb{\abs{x}\to\infty}\bm\phi\sb{c}(x)=0$.
The map $c\mapsto \bm\phi\sb{c}\in H\sp s(\R)$
is $C\sp{\infty }$ for $c\in\varSigma$ and for any $s$.
Consequently,
equation (\ref{gkdv})
admits traveling wave solutions
\begin{equation}
\label{soliton} {\bm u}(x,t)=\bm\phi\sb{c}(x-ct),
\qquad
c\in\varSigma.
\end{equation}
\end{assumption}

In Appendix~\ref{sect-existence} we specify conditions under which
Assumption~\ref{assumption-soliton}
is satisfied.

\begin{figure}[htbp]
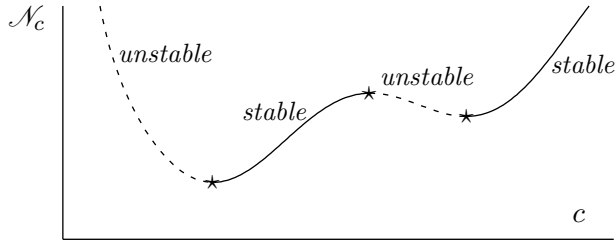

\begin{center}
\input kdvform-n-vs-c.tex
\end{center}
\vskip -30pt
\caption{
Stable and unstable regions
on a possible graph of $\mathscr{N}\sb c$ vs. $c$.
Three critical solitary waves are denoted by stars.}
\label{fig-n-vs-c}
\end{figure}

Let
$\mathscr{N}\sb{c}$ and $I\sb{c}$
denote
$\mathscr{N}(\bm\phi\sb{c})$ and $I(\bm\phi\sb{c})$, respectively.
By Assumption \ref{assumption-soliton},
$\mathscr{N}\sb{c}$ and $I\sb{c}$ are $C\sp\infty$
functions of $c\in\varSigma$.
For the general KdV equation (\ref{gkdv})
with smooth $f({\bm u})$, 
Bona, Suganidis, and Strauss \cite{MR897729} show that
the traveling wave $\bm\phi\sb{c}(x-ct)$
is orbitally stable if
\begin{equation}\label{stability-criterion}
\mathscr{N}\sb{c}'
=\frac{d}{dc}\mathscr{N}\sb{c}
=\frac{d}{dc}\mathscr{N}(\bm\phi\sb{c})
>0
\end{equation}
and unstable if instead $\mathscr{N}\sb{c}'<0$.
See
Figure~\ref{fig-n-vs-c}.
The criterion (\ref{stability-criterion})
coincides with the stability condition
obtained in \cite{MR88g:35169} in the context of
abstract Hamiltonian systems with ${\bf U}(1)$ symmetry
(the theory developed there does not apply
to the generalized Korteweg -- de Vries equation).

\begin{remark}
Note that, as one can readily show,
the amplitude of solitary
waves is monotonically increasing with their speed $c$,
while the momentum $\mathscr{N}\sb c$ does not have to.
\end{remark}

\begin{remark}
For the generalized KdV equations (\ref{kdvk}),
the soliton profiles satisfy the scaling relation
$\bm\phi\sb{c}(x)=c\sp{\frac{1}{p-1}}\bm\phi\sb 1(c\sp{\frac 1 2}x)$.
The values of the momentum
functional that correspond to solitons
with different speeds $c$
are given by
$
\mathscr{N}(\bm\phi\sb{c})
=\mathop{\rm const} c\sp{\frac{2}{p-1}-\frac 1 2}
=\mathop{\rm const} c\sp{\frac{5-p}{2(p-1)}},
$
so that
$
\frac{d}{dc}\mathscr{N}(\bm\phi\sb{c})>0
$
for $p<5$,
in agreement with the stability criterion (\ref{stability-criterion})
derived in \cite{MR897729}.
\end{remark}

In \cite{MR897729} it is stated that   critical traveling waves
$\bm\phi\sb{c\sb\star}(x)$, that is $c\sb\star$ such that
$\mathscr{N}\sb{c\sb\star}'=0$, are unstable as a consequence of the
claim that the set $\{c\sothat \bm\phi\sb{c}{\rm \ is\ stable}
\}$ is open.
This claim however is left unproved in
\cite{MR897729}.
Moreover, this is not true in general.
(This is demonstrated by the dynamical system in $\R^2$
described in the polar coordinates
by $\dot\theta=\sin\theta$, $\dot r=0$.
The set of stationary states is the line $y=0$;
the subset of stable stationary points, $x\le 0$,
is closed.)
The question of stability of critical traveling waves has been left open.
We address this question in this paper,
proving the instability under certain rather generic assumptions.
This result is the analog of \cite{MR1995870}
for the generalized Korteweg -- de Vries equation (\ref{gkdv}).

\begin{remark}
We will not consider the $L\sp 2$-critical KdV equation
given by (\ref{kdvk}) with $p=5$,
when $\mathscr{N}\sb{c}=\const$.
In this case, the solitons
are not only unstable but also exhibit a blow-up behavior.
This blow-up is considered in a series of papers by
Martel and Merle \cite{MR1824989,MR1829643,MR1896235,MR1888800}.
\end{remark}

The analysis of the instability of critical solitary waves
(with no linear instability)
requires better control of the growth
of a particular perturbation.
We achieve this employing the asymptotic stability methods.
Pego and Weinstein \cite{MR1289328}
proved that
the traveling wave solutions to (\ref{kdvk})
for the subcritical values $p=2$, $3$, $4$,
and also $p\in (2,5)\backslash E$ with $E$ a finite and possibly empty set
are asymptotically stable in the weighted spaces.
Their approach was extended in \cite{MR1828318}.
For other deep
results of stability see \cite{ MR1826966,MR2109467}.
The proofs extend,
under certain spectral hypotheses,
to solitary solutions to a generalized KdV equation (\ref{gkdv})
with $c$ such that $\mathscr{N}\sb{c}'>0$.

Substituting ${\bm u}(x,t)=\bm\phi\sb{c}(x-ct)+\bm\rho(x-ct,t)$ into
(\ref{gkdv}) and discarding terms nonlinear in $\bm\rho$,
we get the linearization at $\bm\phi\sb{c}$:
\begin{equation}\label{linearization}
\p\sb{t}\bm\rho
=\p\sb{x}\bigl(-\p\sb{x}\sp{2}\bm\rho
+f'(\bm\phi\sb{c})\bm\rho+c\bm\rho\bigr)
\equiv J\mathcal{H}\sb{c}\bm\rho,
\end{equation}
where
\begin{equation}\label{def-j-h}
J=\p\sb{x},
\qquad
\mathcal{H}\sb{c}=-\p\sb{x}^2+f'(\bm\phi\sb{c})+c.
\end{equation}
In (\ref{linearization}), both $\bm\phi\sb{c}(\cdot)$ and
$\bm\rho(\cdot,t)$ are evaluated at $x-ct$, but we change variable
and write $x$ instead.

\begin{figure}[htbp]
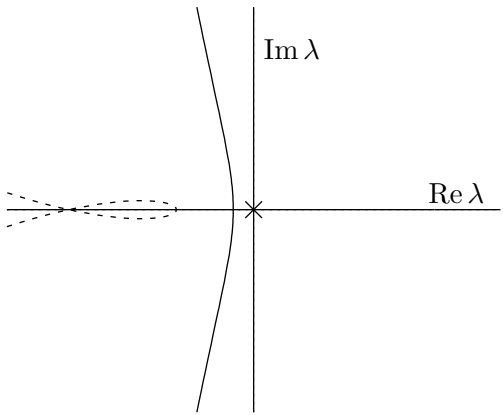

\begin{center}
\input kdvform-spectrum.tex
\end{center}
\vskip -30pt
\caption{Essential spectrum of $J\mathcal{H}\sb{c}$,
$c=1$
in the exponentially weighted space $L\sb\mu\sp{2}(\R)$
for $\mu=0.1<\sqrt{c/3}$ (solid) and $\mu=0.65>\sqrt{c/3}$ (dashed).
}
\label{fig-spectrum}
\end{figure}

The essential spectrum of $J\mathcal{H}\sb{c}$
in $L\sp{2}(\R)$ coincides with the imaginary axis.
$\lambda=0$ is an
eigenvalue
(with $\p\sb x\bm\phi\sb c$ being the corresponding eigenvector).
To use the asymptotic stability methods from \cite{MR1289328},
we will consider the action
of $J\mathcal{H}\sb{c}$ in the exponentially weighted spaces.
For $s\in\R$ and $\mu\ge 0$, we define
\begin{equation}\label{weighted-space}
H\sb\mu\sp{s}(\R)=\left\{\bm\psi\in H\sp{s}\sb{loc}(\R)
\sothat
e\sp{\mu x}\bm\psi(x)\in H\sp{s}(\R)\right\},
\qquad\mu\ge 0,
\end{equation}
where $H\sp s(\R)$
is the standard Sobolev space of order $s$.
We also denote $L\sb\mu\sp{2}(\R)=H\sb\mu\sp{0}(\R).$
We define the operator
$A\sp{\mu}\sb{c}=e^{\mu x}\circ J\mathcal{H}\sb{c}\circ e^{-\mu x}$,
where $e^{\pm\mu x}$ are understood as the operators of multiplication
by the corresponding functions,
so that the action of  $J\mathcal{H}\sb{c}$ in
$L\sb\mu\sp{2}(\R)$ corresponds to the action of
$A\sp{\mu}\sb{c}$ in $L^2(\R)$.
The explicit form of $A\sp{\mu}\sb{c}$ is
\begin{equation}\label{conjugate}
A\sp{\mu}\sb{c}
=e^{\mu x}\circ J\mathcal{H}\sb{c}\circ e^{-\mu x}
=
(\p\sb{x}-\mu)
\big[
-(\p\sb{x}-\mu)^2+c-f'(\bm\phi\sb{c})
\big].
\end{equation}
The domain of $A\sp{\mu}\sb{c}$
is given by $D(A\sp{\mu}\sb{c})=H^3(\R)$.
Since the operator $[\p\sb{x}-\mu]f'(\bm\phi\sb{c})$
is relatively compact with respect to
$\mathscr{A}\sp{\mu}\sb{c}
=-(\p\sb{x}-\mu)^3+ c(\p\sb{x}-\mu)$,
the essential spectrum of $A\sp{\mu}\sb{c}$
coincides with that of $\mathscr{A}\sp{\mu}\sb{c}$
and is given by
\begin{equation}\label{continuous}
\sigma\sb{\rm e}(A\sp{\mu}\sb{c})
=\sigma\sb{\rm e}(\mathscr{A}\sp{\mu}\sb{c})
=\left\{\lambda\in\C
\sothat\lambda=\lambda\sb{\rm cont}(k)=(\mu-i k)^3-c(\mu-i k),
\quad k\in\R
\right\}.
\end{equation}
The essential spectrum of
$\mathscr{A}\sp{\mu}\sb{c}$ is located in the left half-plane for
$0<\mu<\sqrt{c}$
and is simply connected for $0<\mu<\sqrt{c/3}$;
see Figure~\ref{fig-spectrum}.

We need assumptions about the existence and properties
of a critical wave.

\begin{assumption}\label{assumption-existence}
There exists $c\sb\star\in\varSigma\backslash\p\varSigma$,
$c\sb\star>0$,
such that $\mathscr{N}\sb{c\sb\star}'=0$.
\end{assumption}

\begin{remark}\label{remark-critical-waves}
Let us give examples of the nonlinearities
that lead to the existence of critical solitary waves.
Take $f\sb{-}(z)=-Az^p+Bz^q$, with $2<p<q$, $A>0$, $B>0$, 
or
$f\sb{+}(z)=Az^p-Bz^q+Cz^r$, with $2<p<q<r$, $A>0$, $B>0$, $C>0$.
In the case of $f\sb{+}$, we require that
$B$ be sufficiently large so that $f\sb{+}(z)$ takes negative values
on a nonempty interval $I\subset\R\sb{+}$.
Then there will be traveling wave solutions $\bm\phi\sb{c}(x-ct)$
with $c\in(0,c\sb 1)$ (also with $c=0$ in the case of $f\sb{+}$),
for some $c\sb 1>0$.\footnote{The value
of $c\sb 1$ is determined from the system
$f(z\sb 1)+c\sb 1 z\sb 1=0$,
$F(z\sb 1)+c\sb 1 z\sb 1^2/2=0$,
with $F$ the primitive of $f$
such that $F(0)=0$.
See Appendix~\ref{sect-existence} or \cite{MR695535} for more details.}
Elementary computations show that
the value of the momentum $\mathscr{N}\sb{c}$ goes to infinity
as $c\nearrow c\sb 1$.
It also goes to infinity as $c\searrow 0$ if $p>5$
(also if $p=5$ in the case of $f\sb{+}$),
so that there is a global minimum of
$\mathscr{N}\sb{c}$ at some point $c\sb\star\in(0,c\sb 1)$.
\end{remark}

\begin{assumption}\label{assumption-spectrum}
There exists $\mu\sb 0\in(0,\sqrt{c\sb\star}/2)$
such that for $0\le \mu\le\mu\sb 0$
the operator $A\sp{\mu}\sb{c\sb\star}$
has no $L\sp 2$-eigenvalues except $\lambda=0$.
\end{assumption}


\begin{assumption}\label{assumption-non-degeneracy}
At the critical value $c\sb\star$,
the non-degeneracy condition $I\sb{c\sb\star}'\ne 0$ is satisfied.
Here $I\sb{c}=I(\bm\phi\sb c)$
is the value of the mass functional (\ref{mass})
on the traveling wave $\bm\phi\sb c$.
\end{assumption}

\begin{remark}
If
$I\sb{c\sb\star}'=0$,
then the eigenvalue $\lambda=0$ of $J\mathcal{H}\sb{c\sb\star}$
corresponds to a Jordan block larger than $3\times 3$.
We will not consider this situation.
\end{remark}

Our main result is that the critical traveling wave
$\bm\phi\sb{c\sb\star}(x)$
of the generalized KdV equation (\ref{gkdv})
is (nonlinearly) unstable.

\begin{theorem}[Main Theorem]\label{main-theorem}
Let Assumptions~\ref{assumption-soliton},~\ref{assumption-existence},
\ref{assumption-spectrum},
and \ref{assumption-non-degeneracy}
be satisfied,
and that $\bm\phi\sb{c\sb\star}$ is a critical soliton.
Assume that there exists an open neighborhood
$\mathcal{O}(c\sb\star)\subset\varSigma$
of $c\sb\star$
so that 
$\mathscr{N}'\sb{c}$ is strictly negative
and nonincreasing
for $c\in\mathcal{O}(c\sb\star)$, $c>c\sb\star$
(or negative and nondecreasing for $c<c\sb\star$, or both).
Then the critical traveling wave $\bm\phi\sb{c\sb\star}(x)$ is
orbitally unstable.
More precisely,
there exists $\epsilon>0$
such that for any $\delta>0$ 
there exists ${\bm u}\sb 0\in H\sp{1}(\R)$
with $\norm{{\bm u}\sb 0-\bm\phi\sb{c\sb\star}}\sb{H\sp{1}}<\delta$
and $t>0$ so that
\begin{equation}
\inf\sb{s\in\R}
\norm{{\bm u}(\cdot,t)-\bm\phi\sb{c\sb\star}(\cdot-s)}\sb{H\sp{1}}
=\epsilon.
\end{equation}
\end{theorem}

\begin{remark}\label{remark-negative-eta}
For definiteness, we consider the case
when $\mathscr{N}'\sb{c}$ is strictly negative
and nonincreasing for
$c>c\sb\star$, $c\in\mathcal{O}(c\sb\star)$.
The proof for the case
when $\mathscr{N}'\sb{c}$ is strictly negative
and nondecreasing
for $c<c\sb\star$, $c\in\mathcal{O}(c\sb\star)$
is the same.
\end{remark}

Thus, we assume that there exists $\eta\sb 1>0$
such that
\begin{equation}\label{def-eta1}
[c\sb\star,c\sb\star+\eta\sb 1]\subset\varSigma,
\qquad
\mathscr{N}'\sb{c}<0
\quad{\rm for\ }c\in(c\sb\star,c\sb\star+\eta\sb 1]\subset\varSigma.
\end{equation}

\noindent
{\it Strategy of the proof and the structure of the paper.}
\quad
In our proof,
we develop the method of Pego and Weinstein \cite{MR1289328}
and derive the nonlinear bounds
relating the energy estimate and the dissipative estimate
(Lemmas~\ref{lemma-dispersion-0}, \ref{lemma-dispersion}).
We follow a center manifold approach;
that is, we reduce
the infinite-dimensional Hamiltonian system
to a finite dimensional system which contains the main
features of the dynamics.
Specifically, we consider the spectral decomposition near the zero eigenvalue
in Section~\ref{sect-spectral-decomposition}
and a center manifold reduction
is considered in Section~\ref{sect-center-manifold}, this part being
similar to the approach in \cite{MR1995870}.
Estimates in the energy space and in the weighted space for the error terms
are in Section~\ref{sect-estimates} and \ref{sect-nonlinear-estimates}.
In this part of our argument we develop the approach of
\cite{MR1289328}.
In Section~\ref{sect-general-case},
we complete the proof of Theorem~\ref{main-theorem}.
In Section~\ref{sect-normal-forms},
we give an alternative approach to the instability of the
critical traveling wave $\bm\phi\sb{c\sb\star}(x)$
by a normal form argument \cite{MR635782,MR1695170},
under additional hypothesis that the critical point $c\sb\star$
of $\mathscr{N}\sb{c}$ is non-degenerate:
\begin{equation}\label{n-p-p-nonzero}
\mathscr{N}\sb{c\sb\star}''
=\frac{d^2\,\mathscr{N}(\bm\phi\sb{c})}{dc\,^2}\At{c=c\sb\star}
\ne 0.
\end{equation}
The construction of traveling waves
is considered in Appendix~\ref{sect-existence}.
The details on the Fredholm Alternative for $\mathcal{H}\sb{c}$
are in Appendix~\ref{sect-sd}.
An auxiliary technical result is proved in Appendix~\ref{sect-nd}.

\section{Spectral decomposition in $L\sb\mu\sp{2}(\R)$  near $\lambda=0$}
\label{sect-spectral-decomposition}

First, we observe that for any $c\in\varSigma$
(see Assumption~\ref{assumption-soliton}),
the linearization operator
$J\mathcal{H}\sb{c}$ given by (\ref{def-j-h})
satisfies the following relations:
\begin{equation}\label{jh-e1-zero}
\mathcal{H}\sb{c}{\bm e}\sb{1,c}=0,
\qquad
{\rm where}
\quad
{\bm e}\sb{1,c}=-\p\sb{x}\bm\phi\sb{c}(x),
\end{equation}
\begin{equation}\label{jh-e2-e1}
J\mathcal{H}\sb{c}{\bm e}\sb{2,c}={\bm e}\sb{1,c},
\qquad
{\rm where}
\quad
{\bm e}\sb{2,c}=\p\sb{c}\bm\phi\sb{c}(x).
\end{equation}

Let $\mathscr{S}(\R)$ denote the Schwarz space of functions.

\begin{definition}\label{def-Schwarz}
Let $\chi\sb{+}\in C\sp\infty(\R)$
be such that $0\le\chi\sb{+}\le 1$,
$\chi\sb{+}\at{[-1,+\infty)}=0$,
$\chi\sb{+}\at{[0,\infty)}\equiv 1$.
Define
$\mathscr{S}\sb{+,m}(\R)$, $m\ge 0$
to be the set of functions $u\in C\sp\infty(\R)$ such that
$\chi\sb{+} u\in\mathscr{S}(\R)$
and for any $N\in\Z$, $N\ge 0$ there exists
$C\sb N>0$ such that
\[
\abs{u\sp{(N)}(x)}\le C\sb{N}(1+\abs{x})^m.
\]
\end{definition}


Note that for any $m\ge 0$,
$\image(J\mathcal{H}\sb c\at{\mathscr{S}\sb{+,m}(\R)})
\subset\mathscr{S}\sb{+,m}(\R)$.
The algebraic multiplicity of zero eigenvalue
of the operator $J\mathcal{H}\sb c$
considered in $\mathscr{S}\sb{+,m}(\R)$
depends on the values of
$\mathscr{N}\sb{c}'$ and $I\sb{c}'$ as follows.

\begin{proposition}\label{prop-Jordan}
Fix $m\ge 0$,
and consider the operator $J\mathcal{H}\sb c$ in $\mathscr{S}\sb{+,m}(\R)$.

\begin{enumerate}
\item
The eigenvalue $\lambda=0$
is of geometric multiplicity one,
with the kernel generated by ${\bm e}\sb{1,c}$.

\item
Assume that $c\in\varSigma$ is such that
$\mathscr{N}\sb c'\ne 0$.
Then the eigenvalue $\lambda=0$
is of algebraic multiplicity two.

\item
Assume that $c\sb\star\in\varSigma$ is such that
$\mathscr{N}\sb{c\sb\star}'=0$, $I\sb{c\sb\star}'\ne 0$.
Then the eigenvalue $\lambda=0$
is of algebraic multiplicity three.
\end{enumerate}
\end{proposition}

\begin{proof}
First of all we claim that in $\mathscr{S}\sb{+,m}(\R)$
we have $\dim\ker J\mathcal{H}\sb{c}=1$.

The differential equation $\mathcal{H}\sb{c}\bm\psi=0$
has two linearly independent solutions.
According to (\ref{jh-e1-zero}), one of them is
${\bm e}\sb{1,c}$, which is odd and exponentially decaying at
infinity.
The other solution is even and exponentially growing
as $\abs{x}\to\infty$
and hence does not belong to $\mathscr{S}\sb{+,m}(\R)$;
we denote this solution by
$\bm\varXi\sb{c}(x)$.

Observe that if ${\bm v}\in\ker J\mathcal{H}\sb{c}$
then $\mathcal{H}\sb{c}{\bm v}=K$,
${\bm v}\in C\sp\infty(\R)$.
Set ${\bm v}=\frac{K}{c}+{\bm w}$.
Then
$\mathcal{H}\sb{c}{\bm w}=-\frac{K}{c}f'(\bm\phi\sb{c})$.
Since $\left\langle f'(\bm\phi\sb{c}), {\bm e}\sb{1,c}\right\rangle=0$,
by Lemma~\ref{lemma-Fredholm}
there exists a function
${\bm w}\sb 0\in\mathscr{S}\sb{+,m}(\R)$ such that
$\mathcal{H}\sb{c}{\bm w}\sb 0
=
-\frac{K}{c}f'(\bm\phi\sb{c})$.
So ${\bm w}={\bm w}\sb 0+ A\p\sb{x}\bm\phi\sb{c}+B\bm\varXi\sb{c}$,
with $A$ and $B$ constants.
Since
\[
{\bm v}
=\frac{K}{c}+{\bm w}
=\frac{K}{c}
+{\bm w}\sb 0+ A\p\sb{x}\bm\phi\sb{c}+B\bm\varXi\sb{c}
\in\mathscr{S}\sb{+,m}(\R),
\]
we need ${\bm v}(x)\to 0$ for $x\to +\infty$,
and therefore $B=0$ and $K=0$.
Hence, ${\bm v}\in \ker\mathcal{H}\sb c$,
proving that
$\ker J\mathcal{H}\sb c=\ker\mathcal{H}\sb c$.
This proves
Proposition~\ref{prop-Jordan}~({\it i}).

Let us introduce the function
\begin{equation}
\bm\Theta\sb{c}(x)=\int\sb{+\infty}\sp{x}\p\sb{c}\bm\phi\sb{c}(y)\,dy.
\end{equation}
Then
$\p\sb{x}\bm\Theta\sb{c}(x)=\p\sb{c}\bm\phi\sb{c}(x)$,
$\lim\sb{x\to-\infty}\bm\Theta\sb{c}(x)=-I\sb{c}'$,
hence $\Theta\sb{c}\in\mathscr{S}\sb{+,0}(\R)$.
If ${\bm v}$
satisfies
\begin{equation}
J\mathcal{H}\sb{c}{\bm v}=\p\sb{c}\bm\phi\sb{c}(x),
\qquad
\lim\sb{x\to+\infty}{\bm v}(x)=0,
\end{equation}
then ${\bm v}(x)$ is the only solution to the problem
\begin{equation}\label{he-e-infty}
\mathcal{H}\sb{c}{\bm v}=\bm\Theta\sb{c}(x),
\qquad
\lim\sb{x\to+\infty}{\bm v}(x)=0.
\end{equation}
According to Lemma~\ref{lemma-Fredholm}
(see Appendix~\ref{sect-sd}),
if
$
\left\langle{\bm e}\sb{1,c},\bm\Theta\sb{c}\right\rangle
=\left\langle\bm\phi\sb{c},\p\sb{c}\bm\phi\sb{c}\right\rangle
=\mathscr{N}\sb{c}'\ne 0$,
then ${\bm v}(x)$ has exponential growth as $x\to -\infty$:
\begin{equation}\label{eg}
{\bm v}(x)\propto e^{\sqrt{c}\abs{x}},
\qquad
x\to -\infty,
\end{equation}
and therefore does not belong to
$\mathscr{S}\sb{+,m}(\R)$.
This finishes the proof of Proposition~\ref{prop-Jordan}~({\it ii}).

Let us now assume that
$\mathscr{N}\sb{c\sb\star}'=0$
for some $c\sb\star\in\varSigma$.
Then,
again by Lemma~\ref{lemma-Fredholm}
with $m=0$,
there exists
${\bm e}\sb{3,c\sb\star}(x)\in\mathscr{S}\sb{+,0}(\R)$
such that
\begin{equation}\label{inhom}
\mathcal{H}\sb{c\sb\star}{\bm e}\sb{3,c\sb\star}
=\bm\Theta\sb{c\sb\star}(x),
\qquad
\lim\sb{x\to+\infty}{\bm e}\sb{3,c\sb\star}(x)=0.
\end{equation}
Now let us consider ${\bm w}\in C\sp\infty(\R)$
such that
\begin{equation}\label{eg4}
J\mathcal{H}\sb{c\sb\star}{\bm w}={\bm e}\sb{3,c\sb\star},
\qquad
\lim\sb{x\to +\infty}{\bm w}(x)=0.
\end{equation}
Let ${\bm E}(x)=\int\sb{+\infty}\sp{x}{\bm e}\sb{3,c\sb\star}(y)\,dy$;
the function ${\bm w}(x)$ satisfies
$\mathcal{H}\sb{c\sb\star}{\bm w}={\bm E}$.
Taking the pairing of
${\bm E}$ with ${\bm e}\sb{1,c\sb\star}$,
we get:
\begin{eqnarray}
\langle{\bm e}\sb{1,c\sb\star},{\bm E}\rangle
=
-\langle\bm\phi\sb{c\sb\star},{\bm e}\sb{3,c\sb\star}\rangle
=\left\langle
\mathcal{H}\sb{c\sb\star}\p\sb{c}\bm\phi\sb{c\sb\star},{\bm e}\sb{3,c\sb\star}
\right\rangle
=\left\langle
\p\sb{c}\bm\phi\sb{c\sb\star}, \mathcal{H}\sb{c\sb\star}{\bm e}\sb{3,c\sb\star}
\right\rangle
\nonumber
\\
=\left\langle
\p\sb{x}\bm\Theta\sb{c\sb\star},\bm\Theta\sb{c\sb\star}
\right\rangle
=\frac{\bm\Theta\sb{c\sb\star}^2}{2}\biggr|\sb{-\infty}\sp{+\infty}
=-
\lim\sb{x\to-\infty}\frac{\bm\Theta\sb{c\sb\star}^2(x)}{2}
=-\frac{(I\sb{c\sb\star}')^2}{2}<0.
\label{fourth}
\end{eqnarray}
(In the first equality,
the boundary term does not appear
because when $x\to\pm\infty$
the function ${\bm E}(x)$ grows at most
algebraically
while $\bm\phi\sb{c}$ decays exponentially.)
By Lemma~\ref{lemma-Fredholm},
since $\langle{\bm e}\sb{1,c\sb\star},{\bm E}\rangle$ is nonzero,
${\bm w}(x)$ grows exponentially as $x\to -\infty$.
This proves that the algebraic multiplicity
of the eigenvalue $\lambda=0$ is exactly three.
\end{proof}

Now we would like to consider $J \mathcal{H}\sb{c}$
in the weighted space $L\sb\mu\sp 2(\R)$,
$\mu>0$.
This is equivalent to considering
$A\sp{\mu}\sb{c}
=e^{\mu x}\circ J \mathcal{H}\sb{c}\circ e^{-\mu x}$
in $L\sp 2(\R)$.
In what follows, we always require that
\begin{equation}\label{mu-is-such}
0<\mu<\min(\mu\sb 0,\mu\sb 1),
\end{equation}
with $\mu\sb 0$ from Assumption~\ref{assumption-spectrum}
and $\mu\sb 1$ from Lemma~\ref{lemma-nd}.

We define
\begin{equation}
{\bm e}\sp{\mu}\sb{j,c}
=e^{\mu x}{\bm e}\sb{j,c},
\quad j=1,\,2;
\qquad
{\bm e}\sp{\mu}\sb{3,c\sb\star}=e^{\mu x}{\bm e}\sb{3,c\sb\star}.
\end{equation}
  From Proposition~\ref{prop-Jordan},
we obtain the following statement:

\begin{corollary}\label{corollary-Jordan}
\begin{enumerate}
\item
If $\mathscr{N}\sb{c}'\ne 0$,
then
the basis for the generalized kernel of $A\sp{\mu}\sb{c}$
in $L\sp{2}(\R)$
is formed by the generalized eigenvectors
$\{{\bm e}\sp{\mu}\sb{1,c},{\bm e}\sp{\mu}\sb{2,c}\}$.
\item
At $c\sb\star$ where $\mathscr{N}\sb{c\sb\star}'=0$,
the basis for the generalized kernel of $A\sp{\mu}\sb{c\sb\star}$
in $L\sp{2}(\R)$
is formed by the generalized eigenvectors
$\{{\bm e}\sp{\mu}\sb{1,c\sb\star},
{\bm e}\sp{\mu}\sb{2,c\sb\star},{\bm e}\sp{\mu}\sb{3,c\sb\star}\}$.
\end{enumerate}
\end{corollary}

\begin{proof}
As follows from Lemma~\ref{lemma-exp}
in Appendix~\ref{sect-existence},
\begin{equation}
\abs{{\bm e}\sb{1,c}(x)}
\le\const e\sp{-\sqrt{c}\abs{x}},
\qquad
x\in\R.
\end{equation}
Applying Lemma~\ref{lemma-decay}
to (\ref{jh-e2-e1})
(for both $x\ge 0$ and $x\le 0$),
we also see that
\begin{equation}\label{e2-mu}
\abs{{\bm e}\sb{2,c}(x)}
\le\const(1+\abs{x})e\sp{-\sqrt{c}\abs{x}},
\qquad
x\in\R.
\end{equation}
It follows that
${\bm e}\sp\mu\sb{1,c}$,
${\bm e}\sp\mu\sb{2,c}\in L\sp 2(\R)$.

If $\mathscr{N}\sb{c}'\ne 0$,
then by (\ref{eg}) $e^{\mu x}{\bm v}(x)\ne L\sp 2(\R)$.

If $\mathscr{N}\sb{c}'=0$ at $c=c\sb\star$,
then
${\bm e}\sb{3,c\sb\star}\in\mathscr{S}\sb{+,0}(\R)$
(belongs to $\mathscr{S}$ for $x\ge 0$ and remains bounded for $x\le 0$).
Moreover, applying Lemma~\ref{lemma-decay}
to (\ref{inhom}),
we see that
\begin{equation}\label{e3-mu}
\abs{{\bm e}\sb{3,c\sb\star}(x)}
\le\const(1+\abs{x})e\sp{-\sqrt{c\sb\star}\,x},
\qquad
x\ge 0.
\end{equation}
It follows that
${\bm e}\sp\mu\sb{3,c\sb\star}\in L\sp 2(\R)$.
As follows from Proposition~\ref{prop-Jordan},
the function $e^{\mu x}{\bm w}(x)$
in (\ref{eg4}) does not belong to $L\sp 2(\R)$,
so the algebraic multiplicity of $\lambda=0$ is precisely $3$.
\end{proof}

\begin{lemma}\label{lemma-e3}
\begin{enumerate}
\item
Let $c\in(c\sb\star,c\sb\star+\eta\sb 1]$.
Then there exists a simple positive eigenvalue $\lambda\sb{c}$
of $A\sp{\mu}\sb{c}$.
This eigenvalue does not depend on $\mu$.
\item
$\lambda\sb{c}$ is a simple eigenvalue
of the operator $J\mathcal{H}\sb{c}$ considered in $L\sp 2(\R)$.
\item
There exists a $C\sp\infty$ extension of ${\bm e}\sb{3,c\sb\star}$
into an interval $[c\sb\star,c\sb\star+\eta\sb 1]$,
\[
c\mapsto {\bm e}\sb{3,c}\in H\sb\mu\sp\infty(\R),
\qquad c\in[c\sb\star,c\sb\star+\eta\sb 1],
\]
so that the frame
\[
\{
{\bm e}\sp{\mu}\sb{j,c}
=e^{\mu x}{\bm e}\sb{j,c}\in H\sp\infty(\R)
\sothat
j=1,\,2,\,3\},
\qquad
c\in[c\sb\star,c\sb\star+\eta\sb 1]
\]
depends smoothly on $c$ (in $L\sp 2$),
$X\sb c\sp\mu=\mathop{\rm span}
\langle
{\bm e}\sb{1,c}\sp\mu,{\bm e}\sb{2,c}\sp\mu,{\bm e}\sb{3,c}\sp\mu
\rangle$
is the invariant subspace of $A\sp{\mu}\sb{c}$,
and
$A\sp{\mu}\sb{c}\at{X\sb c\sp\mu}$
is represented in the frame
$\{{\bm e}\sb{j,c}\sp\mu\}$
by the following matrix:
\begin{equation}\label{x-d}
A\sp{\mu}\sb{c}\at{
X\sb{c}\sp\mu
}=\left[\begin{array}{ccc}0&1&0
\\0&0&1
\\0&0&\lambda\sb{c}\end{array}
\right],
\end{equation}
where $\lambda\sb{c}$ equals
\begin{equation}\label{def-lambda}
\lambda\sb{c}
=-\frac{\mathscr{N}\sb{c}'}{\langle\bm\phi\sb{c},{\bm e}\sb{3,c}\rangle},
\end{equation}
with
$\langle\bm\phi\sb{c},{\bm e}\sb{3,c}\rangle>0$
for $c\in[c\sb\star,c\sb\star+\eta\sb 1]$.
\end{enumerate}
\end{lemma}

\begin{proof}
Due to the restriction (\ref{mu-is-such}) on $\mu$,
the essential spectrum of $A\sp{\mu}\sb{c}$
for $c\ge c\sb\star$
is given by (\ref{continuous})
and is located strictly to the left of the imaginary axis.
By Assumption~\ref{assumption-spectrum},
the discrete spectrum
of $A\sp{\mu}\sb{c\sb\star}$
consists of the isolated eigenvalue $\lambda=0$,
which is of algebraic multiplicity three
by Corollary~\ref{corollary-Jordan}.
We choose a closed contour $\gamma\subset\rho(A\sp{\mu}\sb{c\sb\star})$
in $\C\sp 1$
so that the interval $[0,\varLambda]$ of the real axis
is strictly inside $\gamma$,
where
\begin{equation}\label{def-varlambda}
\varLambda
=\sup\sb{c\in\varSigma}\,\sup\sb{x\in\R}
\abs{f''(\bm\phi\sb c(x))\bm\phi\sb c'(x)}.
\end{equation}

\begin{remark}\label{remark-lambda}
The value of $\varLambda$
is chosen so that
all pure point eigenvalues of the operator $J\mathcal{H}\sb{c}$,
$c\in\varSigma$,
are bounded by $\varLambda$.
Indeed, if $\bm\psi$ satisfies
$J\mathcal{H}\sb c\bm\psi=\lambda\bm\psi$
with $\lambda\in\R$,
then $\bm\psi\in H\sp\infty(\R)$
and can be assumed real-valued.
Therefore, we have:
\[
\lambda
\langle\bm\psi,\bm\psi\rangle
=
\langle\bm\psi,\p\sb x(-\p\sb x^2+f'(\bm\phi\sb c)+c)\bm\psi\rangle
=-\langle\bm\psi',f'(\bm\phi\sb c)\bm\psi\rangle
=-\langle\bm\psi\bm\psi',f'(\bm\phi\sb c)\rangle
=\frac 1 2
\int\sb{\R}
\bm\psi^2\p\sb x f'(\bm\phi\sb c)\,dx,
\]
so that
$\abs{\lambda}\le\sup\sb{x\in\R}
\abs{f''(\bm\phi\sb c(x))\bm\phi\sb c'(x)}/2$.
\end{remark}

We notice that
for $c$ from an open neighborhood of $c\sb\star$,
$\gamma$ belongs
to the resolvent set
$\rho(A\sp{\mu}\sb{c})$.
Indeed, we have:
\begin{equation}\label{resolvent}
\frac{1}{A\sp{\mu}\sb{c}-z}
=
\frac{1}{A\sp{\mu}\sb{c\sb\star}
-z+(A\sp{\mu}\sb{c}-A\sp{\mu}\sb{c\sb\star})}
=
\frac{1}{(A\sp{\mu}\sb{c\sb\star}-z)}
\frac{1}{(1+(A\sp{\mu}\sb{c\sb\star}-z)^{-1}
(A\sp{\mu}\sb{c}-A\sp{\mu}\sb{c\sb\star}))}.
\end{equation}
Since $A\sp{\mu}\sb{c\sb\star}-z$, $z\in\gamma$,
is invertible in $L\sp 2$ and is smoothing of order three,
while $A\sp{\mu}\sb{c}-A\sp{\mu}\sb{c\sb\star}$ depends continuously on $c$
as a differential operator of order $1$,
the operator
$(A\sp{\mu}\sb{c\sb\star}-z)^{-1}
(A\sp{\mu}\sb{c}-A\sp{\mu}\sb{c\sb\star})$
is bounded by $1/2$ as an operator in $L\sp 2$
for all $z\in\gamma$
and for all $c$
sufficiently close to $c\sb\star$.
We assume that $\eta\sb 1>0$ is small enough so that
\begin{equation}
\gamma\in\rho(A\sp{\mu}\sb{c})
\ \ {\rm for}\ \ c\in[c\sb\star,c\sb\star+\eta\sb 1].
\end{equation}

Integrating (\ref{resolvent}) along $\gamma$,
we get a projection
\begin{equation}\label{def-P-c}
P\sp{\mu}\sb{c}=-\frac{1}{2\pi i}\oint\sb\gamma\frac{dz}{A\sp{\mu}\sb{c}-z},
\qquad
c\in[c\sb\star,c\sb\star+\eta\sb 1].
\end{equation}
Since $\rank P\sp{\mu}\sb{c\sb\star}=3$,
we also have
\[
\rank P\sp{\mu}\sb{c}=3,
\qquad
c\in[c\sb\star,c\sb\star+\eta\sb 1].
\]
The three-dimensional spectral subspace
$\range P\sp{\mu}\sb{c\sb\star}$
corresponds to the eigenvalue $\lambda=0$
that has algebraic multiplicity three.
According to Corollary~\ref{corollary-Jordan},
when $\mathscr{N}\sb{c}'\ne 0$,
$\lambda=0$ is of algebraic multiplicity two,
therefore $X\sb{c}\sp\mu\equiv\range P\sp{\mu}\sb{c}$
splits into a two-dimensional spectral subspace
of $A\sp{\mu}\sb{c}$
corresponding to $\lambda=0$
(it is spanned by $\{{\bm e}\sp{\mu}\sb{1,c},{\bm e}\sp{\mu}\sb{2,c}\}$)
and a one-dimensional subspace 
that corresponds to a nonzero eigenvalue.

For $c\in[c\sb\star,c\sb\star+\eta\sb 1]$,
we define
\begin{equation}\label{def-f-tilde}
\tilde{\bm e}\sp{\mu}\sb{3,c}
=P\sp{\mu}\sb{c}{\bm e}\sp{\mu}\sb{3,c\sb\star},
\qquad
c\in[c\sb\star,c\sb\star+\eta\sb 1].
\end{equation}
Note that
$\tilde{\bm e}\sp{\mu}\sb{3,c}\in L\sp 2(\R)$
since $P\sp{\mu}\sb{c}$ is continuous in $L\sp 2$.
In the frame
$\{{\bm e}\sp{\mu}\sb{1,c},{\bm e}\sp{\mu}\sb{2,c},\tilde{\bm e}\sp{\mu}\sb{3,c}\}$
we can write
\begin{equation}\label{a-f}
A\sp{\mu}\sb{c}
\tilde{\bm e}\sp{\mu}\sb{3,c}
=
a\sb{c}{\bm e}\sp{\mu}\sb{1,c}
+b\sb{c}{\bm e}\sp{\mu}\sb{2,c}
+\lambda\sb{c}
\tilde{\bm e}\sp{\mu}\sb{3,c}.
\end{equation}
Since the frame
$\{{\bm e}\sp{\mu}\sb{1,c},{\bm e}\sp{\mu}\sb{2,c},\tilde{\bm e}\sp{\mu}\sb{3,c}\}$
and also
$A\sp{\mu}\sb{c}\tilde{\bm e}\sp{\mu}\sb{3,c}$
depend smoothly on $c$
(as functions from
$[c\sb\star,c\sb\star+\eta\sb 1]$
to $L\sp 2(\R)$; recall that $f$ is smooth),
the coefficients
$a\sb{c}$, $b\sb{c}$, and $\lambda\sb{c}$ are smooth functions
of $c$ for $c\in[c\sb\star,c\sb\star+\eta\sb 1]$.
It is also important to point out that
$a\sb{c}$, $b\sb{c}$, and $\lambda\sb c$
do not depend on $\mu>0$,
since if the relation (\ref{a-f})
holds for certain values of
$a\sb{c}$, $b\sb{c}$, and $\lambda\sb c$
for a particular value $\mu>0$,
then, by the definition of
$A\sp{\mu}\sb{c}$,
${\bm e}\sp{\mu}\sb{1,c}$, ${\bm e}\sp{\mu}\sb{2,c}$,
and $\tilde{\bm e}\sp{\mu}\sb{3,c}$,
the relation (\ref{a-f}) also holds
for $\mu'$ from an open neighborhood of $\mu$.

According to the construction of
${\bm e}\sb{3,c\sb\star}$
in Proposition~\ref{prop-Jordan},
$a\sb{c\sb\star}=\lambda\sb{c\sb\star}=0$ and
$b\sb{c\sb\star}=1$.
We define
\[
{\bm e}\sp{\mu}\sb{3,c}
=\frac{1}{b\sb{c}+a\sb{c}\lambda\sb{c}}(\tilde{\bm e}\sp{\mu}\sb{3,c}
-a\sb{c}{\bm e}\sp{\mu}\sb{2,c}).
\]
Then ${\bm e}\sp{\mu}\sb{3,c}\in L\sp 2(\R)$
for $c\in[c\sb\star,c\sb\star+\eta\sb 1]$.
We compute:
\begin{equation}\label{e3-relation}
A\sp{\mu}\sb{c}{\bm e}\sp{\mu}\sb{3,c}
={\bm e}\sp{\mu}\sb{2,c}+\lambda\sb{c}{\bm e}\sp{\mu}\sb{3,c}.
\end{equation}
Thus, in the frame $\{{\bm e}\sp{\mu}\sb{j,c}\sothat j=1,\,2,\,3\}$ the
operator $A\sp{\mu}\sb{c}\at{\range P\sp{\mu}\sb{c}}$
has the desired matrix form (\ref{x-d}).
Conjugating by means of $e^{\mu x}$ we get a corresponding
frame $\{{\bm e}\sb{j,c}\sothat j=1,\,2,\,3\}$ in $L\sb\mu^2$, with
${\bm e}\sb{3,c}$ satisfying
\begin{equation}\label{e-3-such}
J\mathcal{H}\sb{c}{\bm e}\sb{3,c}
={\bm e}\sb{2,c}+\lambda\sb{c}{\bm e}\sb{3,c},
\qquad
{\bm e}\sb{3,c}\in L\sb\mu\sp 2(\R).
\end{equation}

For
$c\in[c\sb\star,c\sb\star+\eta\sb 1]$
and
$z\notin\sigma(A\sp{\mu}\sb{c})$,
$R\sp{\mu}\sb{c}(z)=(A\sp{\mu}\sb{c}-z)^{-1}$
is a pseudodifferential operator of order $-3$,
hence
$P\sp{\mu}\sb{c}$
is smoothing of order three in the Sobolev spaces $H\sp s(\R)$.
The bootstrapping argument applied to the relations
${\bm e}\sp{\mu}\sb{j,c}=P\sp{\mu}\sb{c}{\bm e}\sp{\mu}\sb{j,c}$
shows that
${\bm e}\sp{\mu}\sb{j,c}\in H\sp\infty(\R)$.
By definition (\ref{weighted-space}),
this means that
\begin{equation}\label{e-3-smooth}
{\bm e}\sb{j,c}
\in H\sb\mu\sp\infty(\R),
\qquad
j=1,\,2,\,3,
\quad
c\in[c\sb\star,c\sb\star+\eta\sb 1].
\end{equation}
Using (\ref{e-3-such}),
we compute:
\[
0=\langle\mathcal{H}\sb{c}{\bm e}\sb{1,c},{\bm e}\sb{3,c}\rangle
=-\langle\mathcal{H}\sb{c}J\bm\phi\sb{c},{\bm e}\sb{3,c}\rangle
=\langle\bm\phi\sb{c},J\mathcal{H}\sb{c}{\bm e}\sb{3,c}\rangle
=\langle\bm\phi\sb{c},{\bm e}\sb{2,c}\rangle
+\lambda\sb{c}\langle\bm\phi\sb{c},{\bm e}\sb{3,c}\rangle,
\quad c\in[c\sb\star,c\sb\star+\eta\sb 1].
\]
We conclude that
\[
\lambda\sb{c}
=-\frac{\langle\bm\phi\sb{c},{\bm e}\sb{2,c}\rangle}
{\langle\bm\phi\sb{c},{\bm e}\sb{3,c}\rangle},
\qquad
c\in[c\sb\star,c\sb\star+\eta\sb 1],
\]
where
$
\langle\bm\phi\sb{c},{\bm e}\sb{2,c}\rangle
=\langle\bm\phi\sb{c},\p\sb{c}\bm\phi\sb{c}\rangle
=\mathscr{N}\sb{c}'<0.
$
Note that
$\langle\bm\phi\sb{c},{\bm e}\sb{3,c}\rangle>0$
for $c\sb\star<c\le c\sb\star+\eta\sb 1$,
since
$
\langle\bm\phi\sb{c\sb\star},{\bm e}\sb{3,c\sb\star}\rangle>0
$
by (\ref{fourth})
and $\langle\bm\phi\sb{c},{\bm e}\sb{3,c}\rangle$
does not change sign
for $c\sb\star<c\le c\sb\star+\eta\sb 1$
(this follows from the inequality
$\abs{\langle\bm\phi\sb{c},{\bm e}\sb{3,c}\rangle}
>
\abs{\mathscr{N}'\sb{c}}/\varLambda>0$;
see Remark~\ref{remark-lambda}).
%
This finishes the proof of the Lemma.
\end{proof}

\begin{remark}\label{remark-spectrum}
According to Assumption~\ref{assumption-spectrum},
we may assume that
$\eta\sb 1$ is small enough so that
for $c\in[c\sb\star,c\sb\star+\eta\sb 1]$
and $0\le\mu\le\mu\sb 0$
there is no discrete spectrum of $A\sp{\mu}\sb{c}$
except $\lambda=0$ and $\lambda=\lambda\sb{c}$.
It follows that
$P\sb{c}\sp\mu$ is the spectral projector
that corresponds to the discrete spectrum
of $A\sb{c}\sp\mu$.
\end{remark}

\begin{lemma}\label{lemma-is}
If $\lambda\sb c>0$,
then ${\bm e}\sb{3,c}\in H\sp\infty(\R)$.
\end{lemma}

\begin{proof}
By Lemma~\ref{lemma-e3},
$\lambda\sb{c}>0$ is an eigenvalue of $J\mathcal{H}\sb{c}$
considered in $L\sp 2(\R)$.
By
(\ref{jh-e1-zero}), (\ref{jh-e2-e1}), and (\ref{e-3-such}),
\begin{equation}\label{def-psi}
\bm\psi\sb c
={\bm e}\sb{c,1}+\lambda\sb c{\bm e}\sb{c,2}+\lambda\sb c^2{\bm e}\sb{c,3}
\in C\sp\infty(\R)
\end{equation}
satisfies $J\mathcal{H}\sb{c}\bm\psi\sb c=\lambda\sb c\bm\psi\sb c$,
and also
$\lim\sb{x\to+\infty}\bm\psi\sb c(x)=0$.
Thus,
$\bm\psi\sb c$
is an eigenvector
of $J\mathcal{H}\sb{c}\bm\psi\sb c$
that corresponds to $\lambda\sb c$.
Therefore, $\bm\psi\sb c\in H\sp\infty(\R)$.
Since
${\bm e}\sb{c,1}$, ${\bm e}\sb{c,2}\in H\sp 1(\R)$
and $\lambda\sb c\ne 0$, the statement of the lemma
follows from the relation (\ref{def-psi}).
\end{proof}

Let us also introduce
the dual basis that consists of eigenvectors
of the adjoint operator
$\left(J\mathcal{H}\sb{c}\right)\sp\ast
=-\mathcal{H}\sb{c}J
=-\mathcal{H}\sb{c}\p\sb{x}$
which we consider in the weighted space
\begin{equation}
L\sp{2}\sb{-\mu}(\R)
=\left\{\bm\psi\in L\sp{2}\sb{loc}(\R)
\sothat e\sp{-\mu x}\bm\psi(x)\in L\sp{2}(\R) \right\},
\quad\mu > 0.
\end{equation}
For any $c\in\varSigma$,
the generalized kernel of $(J\mathcal{H}\sb{c})\sp\ast$
contains at least two linearly independent vectors:
\begin{equation}
-\mathcal{H}\sb{c}\p\sb{x} {\bm g}\sb{1,c}=0,
\qquad -\mathcal{H}\sb{c}\p\sb{x} {\bm g}\sb{2,c}
={\bm g}\sb{1,c},
\end{equation}
where
\begin{eqnarray}
&&{\bm g}\sb{1,c}(x)
=-\int\sb{-\infty}\sp{x}{\bm e}\sb{1,c}(y,c)\,dy=\bm\phi\sb{c}(x),
\\
&&{\bm g}\sb{2,c}(x)
=\int\sb{-\infty}\sp{x}{\bm e}\sb{2,c}(y,c)\,dy
=\int\sb{-\infty}\sp{x}\p\sb{c}\bm\phi\sb{c}(y)\,dy.
\end{eqnarray}
The lower limit of integration ensures that
$
\lim\sb{x \to -\infty} {\bm g}\sb{2,c}(x)=0,
$
so that
${\bm g}\sb{2,c}\in L\sp{2}\sb{-\mu}(\R)$.

\begin{proposition}
Assume that $c\sb\star\in\varSigma$ is such that
$\mathscr{N}\sb{c\sb\star}'=0$,
$I\sb{c\sb\star}'\ne 0$.
The eigenvalue $\lambda=0$
of the operator $-\mathcal{H}\sb{c\sb\star}\p\sb x$
is of algebraic multiplicity three
in $L\sb{-\mu}\sp 2(\R)$,
and there exists
${\bm g}\sb{3,c\sb\star}\in H\sb{-\mu}\sp\infty(\R)$
such that
\[
-\mathcal{H}\sb{c\sb\star}\p\sb{x}{\bm g}\sb{3,c\sb\star}
={\bm g}\sb{2,c\sb\star}.
\]
\end{proposition}

\begin{proof}
The argument repeats the steps of the proof of Proposition~\ref{prop-Jordan}.
The function ${\bm g}\sb{3,c\sb\star}$ is given by
\begin{equation}
{\bm g}\sb{3,c\sb\star}(x)
=-\int\sb{-\infty}\sp x \tilde{\bm e}\sb{3,c\sb\star}(y) d y,
\end{equation}
where
$
\tilde{\bm e}\sb{3,c\sb\star}(x)
$
satisfies
\begin{equation}
\mathcal{H}\sb{c\sb\star}\tilde{\bm e}\sb{3,c\sb\star}
=\int\sb{-\infty}\sp x {\bm e}\sb{2,c\sb\star}(y)\,dy,
\qquad
\lim\sb{x\to-\infty}\tilde{\bm e}\sb{3,c\sb\star}(x)=0.
\end{equation}
Since
$\int\sb{-\infty}\sp x {\bm e}\sb{2,c\sb\star}(y)\,dy$
remains bounded as $x\to+\infty$,
while
$\langle{\bm g}\sb{2,c\sb\star},\bm\phi\sb{c\sb\star}\rangle=0$,
the function
$\tilde{\bm e}\sb{3,c\sb\star}(x)$
remains bounded as $x\to+\infty$.
This follows from Lemma~\ref{lemma-Fredholm}
of Appendix~\ref{sect-sd}
(after the reflection $x\to-x$).
Therefore,
${\bm g}\sb{3,c\sb\star}(x)$
has a linear growth as $x\to +\infty$;
${\bm g}\sb{3,c\sb\star}\in\mathscr{S}\sb{-,1}(\R)$
(defined similarly to $\mathscr{S}\sb{+,1}$
in Definition~\ref{def-Schwarz}).
\end{proof}

As in Lemma~\ref{lemma-e3},
one can show that there is an extension
of ${\bm g}\sb{3,c\sb\star}$
into an interval $[c\sb\star,c\sb\star+\eta\sb 1]$,
\[
c\mapsto{\bm g}\sb{3,c}\in H\sb{-\mu}\sp\infty(\R),
\qquad
c\in[c\sb\star,c\sb\star+\eta\sb 1],
\]
so that, similarly to (\ref{e-3-such}) and (\ref{e-3-smooth}),
\begin{equation}\label{e3-relation1}
-\mathcal{H}\sb{c}\p\sb{x} {\bm g}\sb{3,c}
={\bm g}\sb{2,c}(x) + \lambda\sb{c}{\bm g}\sb{3,c},
\qquad {\bm g}\sb{3,c}\in H\sp\infty\sb{-\mu}(\R),
\quad c\in[c\sb\star,c\sb\star+\eta\sb 1].
\end{equation}

Using the bases
$\{{\bm e}\sb{j,c}\in H\sb\mu\sp\infty(\R)\sothat j=1,\,2,\,3\}$,
$\{{\bm g}\sb{j,c}\in H\sb{-\mu}\sp\infty(\R)\sothat j=1,\,2,\,3\}$,
we can write the projection operator
$e^{-\mu x}\circ P\sp{\mu}\sb{c}\circ e^{\mu x}$
that corresponds to the discrete spectrum of $J\mathcal{H}\sb{c}$
in the following form:
\begin{equation}\label{pu-pc}
(e^{-\mu x}\circ P\sp{\mu}\sb{c}\circ e^{\mu x})
\bm\psi
=
\sum\sb{j,k=1}\sp{3}
\mathcal{T}\sp{j k}\sb{c}
\langle{\bm g}\sb{k,c},\bm\psi\rangle{\bm e}\sb{j,c},
\end{equation}
with
$\mathcal{T}\sp{j k}\sb{c}$ being the inverse of the matrix
\begin{equation}\label{def-t}
\mathcal{T}\sb{c}
=\{\mathcal{T}\sb{j k,c}\}\sb{1\le j,k\le 3},
\qquad
\mathcal{T}\sb{j k,c}=\langle{\bm g}\sb{j,c},{\bm e}\sb{k,c}\rangle,
\qquad c\in[c\sb\star,c\sb\star+\eta\sb 1],
\quad 1\le j,k\le 3.
\end{equation}
Let us introduce the functions
\begin{equation}\label{alpha-beta-gamma}
\alpha\sb{c}=\langle{\bm g}\sb{1,c},{\bm e}\sb{3,c}\rangle,
\qquad
\beta\sb{c}=\langle{\bm g}\sb{2,c},{\bm e}\sb{3,c}\rangle,
\qquad
\gamma\sb{c}=\langle{\bm g}\sb{3,c},{\bm e}\sb{3,c}\rangle.
\end{equation}
Since
${\bm e}\sb{j,c}\in L\sb\mu\sp 2(\R)$
and
${\bm g}\sb{j,c}\in L\sb{-\mu}\sp 2(\R)$,
$\alpha\sb{c}$, $\beta\sb{c}$,
and $\gamma\sb{c}$
are continuous functions of
$c$ for $c\in[c\sb\star,c\sb\star+\eta\sb 1]$.
Recalling that
$
\langle{\bm g}\sb{2,c},{\bm e}\sb{1,c}\rangle
=\langle{\bm g}\sb{2,c},J\mathcal{H}\sb{c}{\bm e}\sb{2,c}\rangle
=-\langle\mathcal{H}\sb{c}J{\bm g}\sb{2,c},{\bm e}\sb{2,c}\rangle
=\langle{\bm g}\sb{1,c},{\bm e}\sb{2,c}\rangle
=\langle\bm\phi\sb{c},\p\sb{c}\bm\phi\sb{c}\rangle
=\mathscr{N}\sb{c}'$,
$\langle{\bm g}\sb{1,c},{\bm e}\sb{1,c}\rangle
=-\langle\bm\phi\sb{c},\p\sb{x}\bm\phi\sb{c}\rangle=0$,
we may write the matrix $\mathcal{T}$ in the following form:
\begin{equation}\label{def-t-explicit}
\mathcal{T}\sb{c}
=\left[
\begin{array}{ccc} 0&\mathscr{N}\sb{c}'&\alpha\sb{c}
\\
\mathscr{N}\sb{c}'&\frac{1}{2}(I\sb{c}')^2&\beta\sb{c}
\\
\alpha\sb{c}&\beta\sb{c}&\gamma\sb{c}\end{array}\right].
\end{equation}
Note that
$\mathcal{T}\sb{c\sb\star}$ is non-degenerate,
because
$\mathscr{N}\sb{c\sb\star}'=0$ by the choice of $c\sb\star$,
while
$\alpha\sb{c\sb\star}
=\langle{\bm g}\sb{1,c\sb\star},{\bm e}\sb{3,c\sb\star}\rangle
=\langle\bm\phi\sb{c\sb\star},{\bm e}\sb{3,c\sb\star}\rangle
=\frac{1}{2}(I\sb{c\sb\star}')^2>0
$
by (\ref{fourth}).

\section{Center manifold reduction}
\label{sect-center-manifold}

We first discuss the existence of a solution ${\bm u}(t)$
that corresponds to perturbed initial data.
We will rely on the well-posedness results due to T.~Kato.

\begin{lemma}\label{lemma-well}
For any $\mu>0$ and
${\bm u}\sb 0\in H\sp{s}(\R)\cap L\sb{2\mu}\sp 2(\R)$
with
$\norm{{\bm u}\sb 0}\sb{H\sp 1}<2\norm{\bm\phi\sb{c\sb\star}}\sb{H\sp 1}$,
there exists a function
\begin{equation}\label{gwp}
{\bm u}(t)\in C([0,\infty),H\sp{s}(\R)\cap L\sb{2\mu}\sp 2(\R)),
\qquad
{\bm u}(0)={\bm u}\sb 0,
\end{equation}
which solves (\ref{gkdv})
for $0\le t<t\sb 1$,
where $t\sb 1$
is finite or infinite, defined by
\begin{equation}\label{def-t1}
t\sb 1=\sup\{T\in\R\sb{+}\cup\{+\infty\}:\;
\norm{{\bm u}(t)}\sb{H\sp 1}
<2\norm{\bm\phi\sb{c\sb\star}}\sb{H\sp 1}
\ {\rm for}\ t\in(0,T)
\}.
\end{equation}
\end{lemma}

\begin{proof}
According to \cite[Theorem 10.1]{MR759907},
(\ref{gkdv}) is globally well-posed
in $H\sp{s}(\R)\cap L\sb{2\mu}\sp 2(\R)$ for any $s\ge 2$, $\mu>0$
(for the initial data with arbitrarily large norm)
if $f$ satisfies
\begin{equation}\label{kato-gwp}
\lim\sb{\abs{z}\to\infty}\abs{z}^{-4} f'(z)\ge 0.
\end{equation}
We modify the nonlinearity $f(z)$
for $\abs{z}>2\norm{\bm\phi\sb{c\sb\star}}\sb{H\sp 1}$
so that (\ref{kato-gwp}) is satisfied;
Let us call this modified nonlinearity $\tilde{f}(z)$.
Thus, for any
${\bm u}\sb 0\in
H\sp{s}(\R)\cap L\sb{2\mu}\sp 2(\R)$
with
$\norm{{\bm u}\sb 0}\sb{H\sp 1}<2\norm{\bm\phi\sb{c\sb\star}}\sb{H\sp 1}$,
there exists a function
\begin{equation}\label{gwp1}
{\bm u}(t)\in C([0,\infty),H\sp{s}(\R)\cap L\sb{2\mu}\sp 2(\R)),
\qquad
{\bm u}(0)={\bm u}\sb 0,
\end{equation}
that solves the equation
with the modified nonlinearity:
\begin{equation}
\label{gkdv-tilde}
\p\sb{t}{\bm u}
=\p\sb{x}\big(-\p\sb{x}^2 {\bm u}+\tilde{f}({\bm u})\big).
\end{equation}
For $0\le t<t\sb 1$,
with $t\sb 1$ defined by (\ref{def-t1}),
one has
$\norm{{\bm u}(t)}\sb{L\sp\infty}
\le\norm{{\bm u}(t)}\sb{H\sp 1}
<2\norm{\bm\phi\sb{c\sb\star}}\sb{H\sp 1}$.
Therefore,
for $0\le t<t\sb 1$,
${\bm u}(t)$
solves both (\ref{gkdv-tilde}) and (\ref{gkdv})
since $\tilde{f}(z)=f(z)$
for $\abs{z}\le 2\norm{\bm\phi\sb{c\sb\star}}\sb{H\sp 1}$.
\end{proof}

We fix $\mu$ satisfying (\ref{mu-is-such}).
For the initial data
${\bm u}\sb 0\in H\sp 2(\R)\cap L\sb{2\mu}\sp 2(\R)$
with $\norm{{\bm u}\sb 0}\sb{H\sp 1}<2\norm{\bm\phi\sb{c\sb\star}}\sb{H\sp 1}$
there is a function
$\bm u\in C([0,\infty),H\sp 2(\R)\cap L\sb{2\mu}\sp 2(\R))$
that solves (\ref{gkdv}) for $0\le t<t\sb 1$,
with $t\sb 1$ from (\ref{def-t1}).
We will approximate the solution ${\bm u}(x,t)$ by a traveling wave
$\bm\phi\sb{c}$ moving with the variable speed $c=c(t)$.
Thus,
we decompose the
solution ${\bm u}(x,t)$ into the traveling wave $\bm\phi\sb{c}(x)$ and
the perturbation $\bm\rho(x,t)$ as follows:
\begin{equation}\label{def-rho}
{\bm u}(x,t)
=\bm\phi\sb{c(t)}\left(x-\xi(t)-\int\sb 0\sp t c(t')\,dt'\right)
+\bm\rho\left(x-\xi(t)-\int\sb 0\sp t c(t')\,dt',\,t\right).
\end{equation}
The functions
$\xi(t)$ and $c(t)$ are yet to be chosen.

Using (\ref{def-rho}), we rewrite the generalized KdV equation
(\ref{gkdv}) as an equation on $\bm\rho$:
\begin{equation}\label{evolution-rho-0}
\dot{\bm\rho}-J\mathcal{H}\sb{c}\bm\rho
=-\dot\xi{\bm e}\sb{1,c}
-\dot c{\bm e}\sb{2,c}
+\dot\xi\p\sb{x}\bm\rho
+J{\bm N},
\end{equation}
with $\mathcal{H}\sb{c}$  given by (\ref{def-j-h})
and with $J{\bm N}$ given by
\begin{equation}\label{def-n}
J{\bm N}
=\p\sb x
\left[
f(\bm\phi\sb{c}+\bm\rho)-f(\bm\phi\sb{c})-\bm\rho f'(\bm\phi\sb{c})
\right],
\end{equation}
where we changed coordinates,
denoting
$y=x-\xi(t)-{\textstyle\int\sb 0\sp{t}c(t')\,d t'}$
by $x$.
By Proposition~\ref{prop-Jordan}~({\it iii}),
the eigenvalue $\lambda=0$
of operator $J\mathcal{H}\sb{c\sb\star}$ in $L\sb\mu\sp 2(\R)$
has algebraic multiplicity three.
We decompose the perturbation $\bm\rho(x,t)$ as
follows:
\begin{equation}\label{decompos2}
\bm\rho(x,t)=\zeta(t){\bm e}\sb{3,c(t)}(x)+\bm\upsilon(x,t),
\end{equation}
where
${\bm e}\sb{3,c}$ is constructed in Lemma~\ref{lemma-e3}.
Note that the inclusions
$\bm\phi\sb{c}
\in H\sp 2(\R)\cap L\sb{2\mu}\sp 2(\R)
\subset H\sb\mu\sp 1(\R)$
and
${\bm e}\sb{3,c}\in H\sb\mu\sp 1(\R)$
show that $\bm\upsilon(\cdot,t)\in H\sb\mu\sp 1(\R)$.

We would like to choose $\xi(t)$, $c(t)=c\sb\star+\eta(t)$, and $\zeta(t)$
so that
\begin{equation}\label{def-upsilon}
\bm\upsilon(x,t)
={\bm u}
\Big(
x+\xi(t)+\int\sb 0\sp t (c\sb\star+\eta(t'))\,dt',\,t
\Big)
-
\bm\phi\sb{c\sb\star+\eta(t)}(x)
-\zeta(t){\bm e}\sb{3,c\sb\star+\eta(t)}(x)
\end{equation}
represents the part of the perturbation
that corresponds to the continuous spectrum of $J\mathcal{H}\sb{c}$.

\begin{proposition}\label{prop-well}
There exist
$\eta\sb 1>0$, $\zeta\sb 1>0$, and $\delta\sb 1>0$
such that if
$\eta\sb 0$ and $\zeta\sb 0$ satisfy
\begin{equation}\label{init-well}
\abs{\eta\sb 0}<\eta\sb 1,
\qquad
\abs{\zeta\sb 0}<\zeta\sb 1,
\qquad
\norm{
\bm\phi\sb{c\sb\star+\eta\sb 0}+\zeta\sb 0{\bm e}\sb{3,c\sb\star+\eta\sb 0}
-\bm\phi\sb{c\sb\star}
}\sb{H\sp{1}}
<\norm{\bm\phi\sb{c\sb\star}}\sb{H\sp 1},
\end{equation}
then
there is $T\sb 1\in\R\sb{+}\cup\{+\infty\}$
such that:
\begin{enumerate}
\item
There exists
$
{\bm u}\in C([0,\infty),H\sp 2(\R)\cap L\sb{2\mu}\sp 2(\R))
$
so that
\begin{equation}\label{def-u-0}
{\bm u}(0)
=\bm\phi\sb{c\sb\star+\eta\sb 0}+\zeta\sb 0{\bm e}\sb{3,c\sb\star+\eta\sb 0}
\end{equation}
and
${\bm u}(t)$ solves (\ref{gkdv}) for $0\le t<T\sb 1$.
\item
There exist functions
\begin{equation}\label{xi-eta-zeta-0}
\xi,\,\eta,\,\zeta\in C\sp 1([0,\infty)),
\qquad
\xi(0)=0,
\quad
\eta(0)=\eta\sb 0,
\quad
\zeta(0)=\zeta\sb 0,
\end{equation}
such that
the function 
$\bm\upsilon(t)$ defined by (\ref{def-upsilon})
satisfies
\begin{equation}\label{upsilon-rules}
e^{\mu x}\bm\upsilon(x,t)\in\ker P\sp{\mu}\sb{c\sb\star+\eta(t)},
\qquad
0\le t<T\sb 1.
\end{equation}
\item
The following inequalities hold
for $0\le t<T\sb 1$:
\begin{equation}\label{at-t3}
\norm{{\bm u}(t)}\sb{H\sp 1}
<2\norm{\bm\phi\sb{c\sb\star}}\sb{H\sp 1},
\quad
\abs{\eta(t)}<\eta\sb 1,
\quad
\abs{\zeta(t)}<\zeta\sb 1,
\quad
\norm{\bm\upsilon(t)}\sb{H\sb\mu\sp 1}<\delta\sb 1.
\end{equation}
\item
If one can not choose $T\sb 1=\infty$,
then at least one of the inequalities
in (\ref{at-t3}) turns into equality at $t=T\sb 1$.
\end{enumerate}
\end{proposition}


\begin{proof}
Since
${\bm u}\sb 0
=\bm\phi\sb{c\sb\star+\eta\sb 0}+\zeta\sb 0{\bm e}\sb{3,c\sb\star+\eta\sb 0}
\in H\sp 2(\R)\cap L\sb{2\mu}\sp 2(\R)$
and the conditions (\ref{init-well}) are satisfied,
by Lemma~\ref{lemma-well},
there is a function
${\bm u}(t)
\in C([0,\infty),H\sp 2(\R)\cap L\sb{2\mu}\sp 2(\R))$
and $t\sb 1\in\R\sb{+}\cup\{+\infty\}$
such that
${\bm u}(t)$ solves (\ref{gkdv}) for $0\le t<t\sb 1$
and, if $t\sb 1<\infty$, then
$\norm{{\bm u}(t\sb 1)}\sb{H\sp 1}=2\norm{\bm\phi\sb{c\sb\star}}\sb{H\sp 1}$.
We thus need
to construct $\xi(t)$, $\eta(t)$, and $\zeta(t)$
so that
$\bm\upsilon(x,t)$ defined by (\ref{def-upsilon})
satisfies the constraints
\begin{equation}\label{constraints}
\langle{\bm g}\sb{1,c\sb\star+\eta(t)}, \bm\upsilon(t)\rangle
=\langle{\bm g}\sb{2,c\sb\star+\eta(t)}, \bm\upsilon(t)\rangle
=\langle{\bm g}\sb{3,c\sb\star+\eta(t)}, \bm\upsilon(t)\rangle=0.
\end{equation}
Let us note that ${\bm v}(0)=0$
by (\ref{def-upsilon}), (\ref{def-u-0}), and (\ref{xi-eta-zeta-0}).
Since
$J\mathcal{H}\sb{c}{\bm e}\sb{3,c}
=\lambda\sb{c}{\bm e}\sb{3,c}+{\bm e}\sb{2,c}$,
\begin{equation}
\p\sb{t}(\zeta{\bm e}\sb{3,c})
-J\mathcal{H}\sb{c}(\zeta{\bm e}\sb{3,c})
=
\dot\zeta{\bm e}\sb{3,c}
+\dot\eta\zeta\p\sb{c}{\bm e}\sb{3,c}
-\zeta(\lambda\sb{c}{\bm e}\sb{3,c}+{\bm e}\sb{2,c}).
\label{de3-dt}
\end{equation}
Therefore,
(\ref{evolution-rho-0}) can be written as the following equation
on
$\bm\upsilon(t)
=\bm\rho-\zeta{\bm e}\sb{3,c}$:
\begin{equation}\label{dot-upsilon-0}
\dot{\bm\upsilon}-J\mathcal{H}\sb{c}\bm\upsilon
=-
\dot\xi {\bm e}\sb{1,c}-\left(\dot{\eta}-\zeta
\right) {\bm e}\sb{2,c} -(\dot\zeta-\lambda\sb{c}\zeta) {\bm e}\sb{3,c}
-\dot\eta\zeta
\p\sb{c}{\bm e}\sb{3,c}
+\dot\xi\p\sb{x}\bm\rho
+J{\bm N}.
\end{equation}
Differentiating the constraints (\ref{constraints})
and using the evolution equation (\ref{dot-upsilon-0}),
we derive the center manifold reduction:
\begin{equation}\label{cm-reduction1}
\mathcal{T}\sb{c}
\left[
\begin{array}{c}
\dot\xi
\\
\dot{\eta}-\zeta
\\
\dot{\zeta}-\lambda\sb{c}\zeta\end{array}
\right]
-
\dot\eta\left[\begin{array}{cc}
\langle\p\sb{c}{\bm g}\sb{1,c},\bm\upsilon\rangle
\\
\langle\p\sb{c}{\bm g}\sb{2,c},\bm\upsilon\rangle
\\
\langle\p\sb{c}{\bm g}\sb{3,c},\bm\upsilon\rangle
\end{array}\right]
=
-\dot\eta\zeta
\left[\begin{array}{c}
\langle{\bm g}\sb{1,c},\p\sb{c}{\bm e}\sb{3,c}\rangle
\\
\langle{\bm g}\sb{2,c},\p\sb{c}{\bm e}\sb{3,c}\rangle
\\
\langle{\bm g}\sb{3,c},\p\sb{c}{\bm e}\sb{3,c}\rangle
\end{array}\right]
+
\dot\xi
\left[\begin{array}{c}
\langle{\bm g}\sb{1,c},\p\sb{x}\bm\rho\rangle
\\
\langle{\bm g}\sb{2,c},\p\sb{x}\bm\rho\rangle
\\
\langle{\bm g}\sb{3,c},\p\sb{x}\bm\rho\rangle
\end{array}\right]
+
\left[\begin{array}{c}
\langle{\bm g}\sb{1,c},J{\bm N}\rangle
\\
\langle{\bm g}\sb{2,c},J{\bm N}\rangle
\\
\langle{\bm g}\sb{3,c},J{\bm N}\rangle
\end{array}\right],
\end{equation}
where the matrix
$\mathcal{T}\sb{c}$
is given by (\ref{def-t}).
The above can be rewritten as
\begin{equation}\label{cm-reduction11}
\mathcal{S}
\left[\begin{array}{cc}\dot\xi
\\
\dot{\eta}-\zeta
\\
\dot{\zeta}-\lambda\sb{c}\zeta \end{array}\right]
=
\left[\begin{array}{cc}
-\zeta^2\langle{\bm g}\sb{1,c},\p\sb{c}{\bm e}\sb{3,c}\rangle
+\zeta\langle\p\sb{c}{\bm g}\sb{1,c},\bm\upsilon\rangle
+\langle{\bm g}\sb{1,c},J{\bm N}\rangle
\\
-\zeta^2\langle{\bm g}\sb{2,c},\p\sb{c}{\bm e}\sb{3,c}\rangle
+\zeta\langle\p\sb{c}{\bm g}\sb{2,c},\bm\upsilon\rangle
+\langle{\bm g}\sb{2,c},J{\bm N}\rangle
\\
-\zeta^2\langle{\bm g}\sb{3,c},\p\sb{c}{\bm e}\sb{3,c}\rangle
+\zeta\langle\p\sb{c}{\bm g}\sb{3,c},\bm\upsilon\rangle
+\langle{\bm g}\sb{3,c},J{\bm N}\rangle
\end{array}\right],
\end{equation}
where $c=c\sb\star+\eta$ and
\begin{equation}\label{def-s}
\mathcal{S}(\eta,\zeta,\bm\upsilon)
=
\mathcal{T}\sb{c}
+
\left[\begin{array}{ccc}
-\langle{\bm g}\sb{1,c},\p\sb{x}(\zeta{\bm e}\sb{3,c}+\bm\upsilon)\rangle
&\zeta\langle{\bm g}\sb{1,c},\p\sb{c}{\bm e}\sb{3,c}\rangle
-\langle\p\sb{c}{\bm g}\sb{1,c},\bm\upsilon\rangle&0
\\
-\langle{\bm g}\sb{2,c},\p\sb{x}(\zeta{\bm e}\sb{3,c}+\bm\upsilon)\rangle&
\quad
\zeta\langle
{\bm g}\sb{2,c},\p\sb{c}{\bm e}\sb{3,c}\rangle
-\langle\p\sb{c}{\bm g}\sb{2,c},\bm\upsilon\rangle
\quad&0
\\
-\langle{\bm g}\sb{3,c},\p\sb{x}(\zeta{\bm e}\sb{3,c}+\bm\upsilon)\rangle&
\zeta\langle
{\bm g}\sb{3,c},\p\sb{c}{\bm e}\sb{3,c}\rangle
-\langle\p\sb{c}{\bm g}\sb{3,c},\bm\upsilon\rangle
&0
\end{array}\right].
\end{equation}
Note that
the matrix $\mathcal{S}(\eta,\zeta,\bm\upsilon)$
depends continuously on
$(\eta,\zeta,\bm\upsilon)\in\R^2\times H\sb{\mu}\sp 1(\R)$.
Since the matrix $\mathcal{T}\sb{c\sb\star}$
is non-singular
(see (\ref{def-t-explicit})),
the matrix $\mathcal{S}(\eta,\zeta,\bm\upsilon)$ is invertible
for sufficiently small values of
$\abs{\eta}$,
$\abs{\zeta}$, and $\norm{\bm\upsilon}\sb{H\sb\mu\sp 1}$.

Thus,
there exist
$\eta\sb 1>0$, $\zeta\sb 1>0$,
and $\delta\sb 1>0$
so that the matrix $\mathcal{S}(\eta,\zeta,\bm\upsilon)$
is invertible if
\begin{equation}\label{delta-0}
\abs{\eta}\le 2\eta\sb 1,
\qquad
\abs{\zeta}\le 2\zeta\sb 1,
\qquad
\norm{\bm\upsilon}\sb{H\sb\mu\sp 1}\le 2\delta\sb 1.
\end{equation}
For such
$\eta$, $\zeta$, and $\bm\upsilon$,
we can write
\begin{equation}\label{xi-eta-zeta-system}
\left[
\begin{array}{c}
\dot\xi
\\
\dot\eta-\zeta
\\
\dot\zeta-\lambda\sb{c}\zeta
\end{array}
\right]
=
\left[
\begin{array}{c}
R\sb 0(\eta,\zeta,\bm\upsilon)
\\
R\sb 1(\eta,\zeta,\bm\upsilon)
\\
R\sb 2(\eta,\zeta,\bm\upsilon)
\end{array}
\right],
\end{equation}
where the right-hand-side is given by
\begin{equation}\label{def-r}
\left[
\begin{array}{c}
R\sb 0(\eta,\zeta,\bm\upsilon)
\\
R\sb 1(\eta,\zeta,\bm\upsilon)
\\
R\sb 2(\eta,\zeta,\bm\upsilon)
\end{array}
\right]
=\mathcal{S}(\eta,\zeta,\bm\upsilon)\sp{-1}
\left[\begin{array}{cc}
-\zeta^2\langle{\bm g}\sb{1,c},\p\sb{c}{\bm e}\sb{3,c}\rangle
+\zeta\langle\p\sb{c}{\bm g}\sb{1,c},\bm\upsilon\rangle
+\langle{\bm g}\sb{1,c},J{\bm N}\rangle
\\
-\zeta^2\langle{\bm g}\sb{2,c},\p\sb{c}{\bm e}\sb{3,c}\rangle
+\zeta\langle\p\sb{c}{\bm g}\sb{2,c},\bm\upsilon\rangle
+\langle{\bm g}\sb{2,c},J{\bm N}\rangle
\\
-\zeta^2\langle{\bm g}\sb{3,c},\p\sb{c}{\bm e}\sb{3,c}\rangle
+\zeta\langle\p\sb{c}{\bm g}\sb{3,c},\bm\upsilon\rangle
+\langle{\bm g}\sb{3,c},J{\bm N}\rangle
\end{array}\right].
\end{equation}
Assume that $\eta\sb 0$ and $\zeta\sb 0$
are such that the conditions (\ref{init-well}) are satisfied.
Let $\varrho\sb 0\in C\sb{comp}\sp\infty(\R)$
be such that
$0\le\varrho\sb 0(s)\le 1$,
$\varrho\sb 0(s)\equiv 1$ for $\abs{s}\le 1$,
and
$\varrho\sb 0(s)\equiv 0$ for $\abs{s}\ge 2$.
Define a continuous matrix-valued function
$\tilde{\mathcal{S}}:\;\R^2\times H\sb\mu\sp 1\to GL(3)$
by
\[
\tilde{\mathcal{S}}(\eta,\zeta,\bm\upsilon)
=\mathcal{S}
(\varrho\eta,
\varrho\zeta,
\varrho\bm\upsilon),
\qquad
{\rm where}
\quad
\varrho=
\varrho\sb 0(\fra{\eta}{\eta\sb 1})
\varrho\sb 0(\fra{\zeta}{\zeta\sb 1})
\varrho\sb 0(\fra{\norm{\bm\upsilon}\sb{H\sb\mu\sp 1}}{\delta\sb 1}).
\]
This function coincides with $\mathcal{S}$
(defined in (\ref{def-s}))
for $\abs{\eta}<\eta\sb 1$,
$\abs{\zeta}<\zeta\sb 1$,
and $\norm{\bm\upsilon}\sb{H\sb\mu\sp 1}<\delta\sb 1$,
and has uniformly bounded inverse.
The system (\ref{xi-eta-zeta-system})
with the right-hand side 
as in (\ref{def-r})
but with $\tilde{\mathcal{S}}$ instead of $\mathcal{S}$,
and with $\bm\upsilon$ given by the ansatz (\ref{def-upsilon}),
defines differentiable functions $\xi(t)$, $\eta(t)$, and $\zeta(t)$
for all $t\ge 0$.
Note that $\bm\upsilon(t)$
defined by (\ref{def-upsilon})
is a continuous function of time,
and is valued in $H\sb\mu\sp 1(\R)$
since so are ${\bm u}$, $\bm\phi\sb c$, and ${\bm e}\sb{3,c}$.
Define $t\sb 2\in\R\sb{+}\cup\{+\infty\}$ by
\begin{equation}\label{def-t2}
t\sb 2
=\sup\{T\in\R\sb{+}\cup\{+\infty\}
\sothat
\abs{\eta(t)}<\eta\sb 1,
\ \abs{\zeta(t)}<\zeta\sb 1,
\ \norm{\bm\upsilon(\cdot,t)}\sb{H\sb\mu\sp 1}<\delta\sb 1
\ {\rm for}\ t\in(0,T)
\}.
\end{equation}
For $t\in(0,t\sb 2)$,
the solution
$(\xi(t),\eta(t),\zeta(t))$
also solves (\ref{xi-eta-zeta-system}),
since the inequalities
$
\abs{\eta(t)}<\eta\sb 1,
$
$
\abs{\zeta(t)}<\zeta\sb 1,
$
and
$
\norm{\bm\upsilon(\cdot,t)}\sb{H\sb\mu\sp 1}<\delta\sb 1
$
ensure that $\tilde{\mathcal{S}}$ coincides with $\mathcal{S}$.
Thus, Proposition~\ref{prop-well}
is proved with
\begin{equation}\label{def-time}
T\sb 1=\min(t\sb 1,\,t\sb 2)\in\R\sb{+}\cup\{+\infty\},
\end{equation}
where $t\sb 1$, $t\sb 2$ are from (\ref{def-t1}) and (\ref{def-t2}).
\end{proof}

\section{Energy and dissipative estimates}
\label{sect-estimates}

We will adapt the analysis from \cite{MR1289328}.
In this section,
we formulate two Lemmas
that are the analog of \cite[Proposition 6.1]{MR1289328}.
Lemma~\ref{lemma-energy}
is based on the energy conservation
and allows to control
$\norm{\bm\rho}\sb{H\sp 1}$
in terms of $\norm{\bm\upsilon}\sb{H\sb\mu\sp{1}}$.
Lemma~\ref{lemma-dispersion}
bounds $\norm{\bm\upsilon}\sb{H\sb\mu\sp{1}}$
in terms of $\norm{\bm\rho}\sb{H\sp 1}$ and is based on dissipative
estimates on the semigroup generated by $A\sp{\mu}\sb{c}$
(see Lemma~\ref{lemma-dispersion-0}).

Let $\eta\sb 1>0$, $\zeta\sb 1>0$, and $\delta\sb 1>0$
be not larger than in Proposition~\ref{prop-well},
and assume that $\delta\sb 1$
satisfies
\begin{equation}\label{delta-small-1}
\delta\sb 1
<\frac{\min(1,c\sb\star)}
{4\sup\sb{\abs{z}\le 2\norm{\bm\phi\sb{c\sb\star}}\sb{H\sp 1}}\abs{f''(z)}}.
\end{equation}
Let 
$\eta\sb 0>0$ and $\zeta\sb 0$
be such that
the conditions (\ref{init-well}) are satisfied.
According to Proposition~\ref{prop-well},
there exists $T\sb 1\in\R\sb{+}\cup\{+\infty\}$ such that
there is a solution
${\bm u}\in C((0,T\sb 1),H\sp{2}(\R)\cap L\sb{2\mu}\sp 2(\R))$
to (\ref{gkdv})
with the initial data
${\bm u}\sb 0=\bm\phi\sb{c\sb\star+\eta\sb 0}
+\zeta\sb 0{\bm e}\sb{3,c\sb\star+\eta\sb 0}$,
and functions $\xi(t)$, $\eta(t)$, and $\zeta(t)$
and $\bm\upsilon(t)$ (given by (\ref{def-upsilon})),
defined for $0\le t<T\sb 1$,
such that (\ref{upsilon-rules}) and (\ref{at-t3})
are satisfied.
For given $\eta\sb 0$ and $\zeta\sb 0$,
define the following function of $\eta$:
\begin{equation}\label{def-y}
\mathscr{Y}(\eta)
=
\norm{\bm\rho\sb 0}\sb{H\sp{1}}
+\norm{\bm\rho\sb 0}\sb{H\sp{1}}\sp{1/2}\abs{\eta-\eta\sb 0}^{1/2}
+\abs{\mathscr{N}\sb{c\sb\star+\eta}
-\mathscr{N}\sb{c\sb\star+\eta\sb 0}}\sp{1/2},
\qquad
\bm\rho\sb 0
\equiv\zeta\sb 0{\bm e}\sb{3,c\sb\star+\eta\sb 0}.
\end{equation}

\begin{lemma}\label{lemma-energy}
There exists
$C\sb 1>0$
such that if
at some moment $0\le t<T\sb 1$
\[
\norm{\bm\rho(t)}\sb{H\sp{1}} \le\delta\sb 1,
\]
then
\begin{equation}\label{bound-rho0}
\norm{\bm\rho(t)}\sb{H\sp{1}}
\le
C\sb 1\left(
\mathscr{Y}(\eta(t))
+\abs{\zeta(t)}+\norm{\bm\upsilon(t)}\sb{H\sb\mu\sp{1}}
\right),
\end{equation}
where $\mathscr{Y}(\eta)$ is given by (\ref{def-y}).
\end{lemma}

\begin{proof}
Let us introduce the effective Hamiltonian $\mathscr{L}\sb{c}$:
\begin{equation}\label{def-lc}
\mathscr{L}\sb{c}({\bm u})=E({\bm u})+c \mathscr{N}({\bm u}),
\qquad
\mathscr{L}\sb{c}'(\bm\phi\sb{c})
=E'(\bm\phi\sb{c})+c\mathscr{N}'(\bm\phi\sb{c})=0,
\qquad
\mathscr{L}\sb{c}''(\bm\phi\sb{c})=\mathcal{H}\sb{c},
\end{equation}
where $E$ and $\mathscr{N}$ are the energy and momentum
functionals defined in
(\ref{hamiltonian}) and (\ref{momentum}).
Using the Taylor series expansion for $\mathscr{L}\sb{c}$
at $\bm\phi\sb{c}$, we have:
\begin{eqnarray}\label{with-error-term}
&&\mathscr{L}\sb{c}({\bm u}(t))
=\mathscr{L}\sb{c}(\bm\phi\sb{c})
+\frac{1}{2}\left\langle\bm\rho,\mathcal{H}\sb{c}\bm\rho\right\rangle
+\int\sb{\R} g(\bm\phi\sb{c},\bm\rho) \bm\rho^3\,dx
\nonumber
\\
&&
=\mathscr{L}\sb{c}(\bm\phi\sb{c})
+\frac{1}{2}\langle\bm\rho,(-\p\sb{x}^2+c)\bm\rho\rangle
+\frac{1}{2}\langle\bm\rho,f'(\bm\phi\sb{c})\bm\rho\rangle
+\int\sb{\R} g(\bm\phi\sb{c},\bm\rho) \bm\rho^3\,dx,
\end{eqnarray}
where
\begin{equation}\label{def-gg}
g(\bm\phi\sb{c},\bm\rho)
=\frac{1}{2}\int\sb 0\sp{1}(1-s)^2 f''(\bm\phi\sb{c}+s\bm\rho)\,ds.
\end{equation}
For the second term in (\ref{with-error-term}),
there is the following bound from below:
\begin{equation}\label{quadratic-part-p}
\frac 1 2\int\sb{\R} \left((\p\sb{x}\bm\rho)^2+c\bm\rho^2
\right)\,dx
\ge
m\norm{\bm\rho}\sb{H\sp{1}}^2,
\qquad
m=\frac{1}{2}\min(1,c\sb\star)>0.
\end{equation}
There is the following bound for the third term
in the right-hand side of (\ref{with-error-term}):
\begin{equation}
\frac 1 2
\int\sb{\R}\abs{f'(\bm\phi\sb{c})}\bm\rho^2\,dx
\le
\frac 1 2
\norm{e\sp{-2\mu x} f'(\bm\phi\sb{c})}\sb{L\sp\infty}
\norm{\bm\rho}\sb{L\sb\mu\sp{2}}^2
\le
\frac{b}{2}
\left[
\abs{\zeta}\norm{{\bm e}\sb{3,c}}\sb{L\sb\mu\sp 2}
+\norm{\bm\upsilon(t)}\sb{L\sb\mu\sp{2}}\right]^2,
\label{thanks-weights-p}
\end{equation}
where
$b=\sup\sb{c\in[c\sb\star,c\sb\star+\eta\sb 1]}
\norm{e^{-2\mu x}f'(\bm\phi\sb c)}\sb{L\sp\infty}
<\infty$
due to (\ref{mu-is-such}),
the assumption
(\ref{f-p-zero})
that $f'(0)=0$,
and due to Lemma~\ref{lemma-exp} from Appendix~\ref{sect-existence}.
We bound the last term in (\ref{with-error-term})
by
\begin{equation}
\int\sb{\R}\abs{g(\bm\phi\sb{c},\bm\rho) \bm\rho^3}\,dx
\le 
\norm{g(\bm\phi\sb c,\bm\rho)}\sb{L\sp\infty}
\norm{\bm\rho}\sb{H\sp{1}}^3
\le
\delta\sb 1
\norm{g(\bm\phi\sb c,\bm\rho)}\sb{L\sp\infty}
\norm{\bm\rho}\sb{H\sp{1}}^2.
\end{equation}
According to (\ref{delta-small-1}),
$g$ from (\ref{def-gg}) satisfies
$
\delta\sb 1\norm{g(\bm\phi\sb c,\bm\rho)}\sb{L\sp\infty}
<\frac{\min(1,c\sb\star)}{4}
=\frac{m}{2}
$,
and this leads to
\begin{equation}\label{bound-error-term}
\int\sb{\R}\abs{g(\bm\phi\sb{c},\bm\rho) \bm\rho^3}\,dx
\le\frac{m}{2}\norm{\bm\rho}\sb{H\sp{1}}^2.
\end{equation}
Combining (\ref{with-error-term})
with the bounds (\ref{quadratic-part-p}), (\ref{thanks-weights-p}),
and (\ref{bound-error-term}),
we obtain:
\[
\frac{m}{2}
\norm{\bm\rho}\sb{H\sp{1}}^2
\le
\abs{\mathscr{L}\sb{c}({\bm u})-\mathscr{L}\sb{c}(\bm\phi\sb{c})}
+
\frac{b}{2}
\left[
\abs{\zeta}
\norm{{\bm e}\sb{3,c}}\sb{L\sb\mu\sp{2}}
+
\norm{\bm\upsilon}\sb{H\sb\mu\sp{1}}
\right]^2,
\]
so that, for some $C>0$,
\begin{equation}\label{rho-h1}
\norm{\bm\rho}\sb{H\sp{1}}
\le
C\left[
\abs{\mathscr{L}\sb{c}({\bm u})-\mathscr{L}\sb{c}(\bm\phi\sb{c})}\sp{1/2}
+\abs{\zeta}+\norm{\bm\upsilon}\sb{H\sb\mu\sp{1}}
\right].
\end{equation}

Now let us estimate
$\abs{\mathscr{L}\sb{c}({\bm u}(t))-\mathscr{L}\sb{c}(\bm\phi\sb{c})}$.
Note that
$\mathscr{L}\sb{c}({\bm u}(t))=\mathscr{L}\sb{c}({\bm u}(0))$
since the value of the energy functional $E$ given by (\ref{hamiltonian})
and the value of the momentum functional $\mathscr{N}$
given by (\ref{momentum})
are conserved along the trajectories of equation (\ref{gkdv}).
Thus, we can write:
\begin{equation}
\label{l-l}
\abs{\mathscr{L}\sb{c}({\bm u}(t))-\mathscr{L}\sb{c}(\bm\phi\sb{c})}
\le\abs{\mathscr{L}\sb{c}({\bm u}(0))
-\mathscr{L}\sb{c}(\bm\phi\sb{c\sb 0})}
+\abs{\mathscr{L}\sb{c}(\bm\phi\sb{c})-\mathscr{L}\sb{c}(\bm\phi\sb{c\sb 0})}.
\end{equation}
Using the definition (\ref{def-lc}) of the functional
$\mathscr{L}\sb{c}$,
we express the first term in the right-hand side of (\ref{l-l}) as
\begin{equation}\label{lc-lc}
\mathscr{L}\sb{c}({\bm u}(0))-\mathscr{L}\sb{c}(\bm\phi\sb{c\sb 0})
=
\mathscr{L}\sb{c\sb 0}({\bm u}(0))-\mathscr{L}\sb{c\sb 0}(\bm\phi\sb{c\sb 0})
+(\eta-\eta\sb 0)(\mathscr{N}({\bm u}(0))-\mathscr{N}(\bm\phi\sb{c\sb 0})).
\end{equation}
Since $\mathscr{L}\sb{c\sb 0}'(\bm\phi\sb{c\sb 0})=0$,
there exists $k>0$ such that
$\abs{
\mathscr{L}\sb{c\sb 0}({\bm u}(0))-\mathscr{L}\sb{c\sb 0}(\bm\phi\sb{c\sb 0})}
\le k \norm{\bm\rho\sb 0}\sb{H\sp{1}}^2$,
where $\bm\rho\sb 0={\bm u}(0)-\bm\phi\sb{c\sb 0}$;
this allows to bound (\ref{lc-lc}) by
\begin{equation}
\label{l-l-first-term}
\abs{\mathscr{L}\sb{c}({\bm u}(0))-\mathscr{L}\sb{c}(\bm\phi\sb{c\sb 0})}
\le\const
(\norm{\bm\rho\sb 0}\sb{H\sp{1}}^2
+\abs{\eta-\eta\sb 0}\norm{\bm\rho\sb 0}\sb{H\sp{1}}).
\end{equation}
For the second term in the right-hand side of (\ref{l-l}),
we have:
\[
\abs{\mathscr{L}\sb{c}(\bm\phi\sb{c})-\mathscr{L}\sb{c}(\bm\phi\sb{c\sb 0})}
\le\abs{E\sb{c}-E\sb{c\sb 0}}+c\abs{\mathscr{N}\sb{c}-\mathscr{N}\sb{c\sb 0}}.
\]
  From the relation
\[
\frac{d}{dc}E\sb{c}=-c \frac{d}{dc}\mathscr{N}\sb{c}
\]
we conclude that
$\abs{E\sb{c}-E\sb{c\sb 0}}
\le\max(c,c\sb 0)\abs{\mathscr{N}\sb{c}-\mathscr{N}\sb{c\sb 0}}$,
since $\mathscr{N}'\sb{c}$
is sign-definite
for $c\sb\star<c\le c\sb\star+\eta\sb 1$
by (\ref{def-eta1}).
Therefore,
there is the following bound for the second term
in the right-hand side of (\ref{l-l}):
\begin{equation}
\label{l-l-second-term}
\abs{\mathscr{L}\sb{c}(\bm\phi\sb{c})-\mathscr{L}\sb{c}(\bm\phi\sb{c\sb 0})}
\le 2\max(c,c\sb 0)\abs{\mathscr{N}\sb{c}-\mathscr{N}\sb{c\sb 0}}.
\end{equation}
Using the bounds (\ref{l-l-first-term})
and (\ref{l-l-second-term})
in the inequality (\ref{l-l}), we obtain:
\[
\abs{\mathscr{L}\sb{c}({\bm u}(t))
-\mathscr{L}\sb{c}(\bm\phi\sb{c})}
\le
\const
\Big(
\norm{\bm\rho\sb 0}\sb{H\sp 1}^2
+\abs{\eta-\eta\sb 0}\norm{\bm\rho\sb 0}\sb{H\sp 1}
+\abs{\mathscr{N}\sb{c}-\mathscr{N}\sb{c\sb 0}}
\Big).
\]
Substituting this result into (\ref{rho-h1}),
we obtain the bound (\ref{bound-rho0}).
\end{proof}

\begin{lemma}[\cite{MR1289328}]\label{lemma-dispersion-0}
Let Assumption~\ref{assumption-spectrum} be satisfied, and pick
$\mu\in(0,\sqrt{c/3})$.
Let $Q\sp{\mu}\sb{c}=I-P\sp{\mu}\sb{c}$,
where $P\sp{\mu}\sb{c}$ introduced in (\ref{def-P-c})
is the spectral projection
that corresponds to the discrete
spectrum of $A\sp{\mu}\sb{c}$
(see Remark~\ref{remark-spectrum}).
Then $A\sp{\mu}\sb{c}$
is the generator of a strongly continuous linear semigroup on $H^s
(\R)$ for any real $s$, and there exist constants $a > 0$ and $b >
0$ such that for all $\bm\upsilon\in L\sp{2}(\R)$ and $t>0$ the
following estimate is satisfied:
\begin{equation}\label{smoothing}
\norm{e\sp{A\sp{\mu}\sb{c}t} Q\sp{\mu}\sb{c}\bm\upsilon}\sb{H\sp{1}}
\le
a t\sp{-1/2} e\sp{-b t}\norm{\bm\upsilon}\sb{L\sp{2}}.
\end{equation}
\end{lemma}

We require that $\eta\sb 1$ be small enough,
so that
\begin{equation}\label{eta-small-q}
\eta\sb 1
\sup\sb{c\in[c\sb\star,c\sb\star+\eta\sb 1]}
\norm{\p\sb{c}Q\sp{\mu}\sb{c}}\sb{H\sp 1\to H\sp 1}
\le \frac 1 2.
\end{equation}

\begin{lemma}\label{lemma-dispersion}
There exists $C\sb 2>0$
such that if
\begin{equation}\label{eta-zeta-delta-small}
\eta\sb 1+\zeta\sb 1+\delta\sb 1<C\sb 2
\end{equation}
and
\begin{equation}\label{uniform-bounds-assumed}
\sup\sb{s\in[0,t]}\abs{\eta(s)}\le\eta\sb 1,
\quad
\sup\sb{s\in[0,t]}\abs{\zeta(s)}\le\zeta\sb 1,
\quad
\sup\sb{s\in[0,t]}
\norm{\bm\rho(s)}\sb{H\sp{1}}
\le\delta\sb 1,
\quad
\sup\sb{s\in[0,t]}
\norm{\bm\upsilon(s)}\sb{H\sb\mu\sp{1}}\le\delta\sb 1,
\end{equation}
then
\begin{equation}\label{bound-v-mu0}
\norm{\bm\upsilon(t)}\sb{H\sb\mu\sp{1}}
\le C\sb 2
\sup\sb{s\in[0,t]}
\left[\zeta^2(s)+\abs{\zeta(s)}\norm{\bm\rho(s)}\sb{H\sp{1}}
\right].
\end{equation}
\end{lemma}

\begin{proof}
Using the center manifold
reduction (\ref{xi-eta-zeta-system}),
we rewrite the evolution equation (\ref{dot-upsilon-0})
in the following form:
\begin{equation}
\label{evolution-upsilon}
\dot{\bm\upsilon}-J\mathcal{H}\sb{c}\bm\upsilon
=-R\sb 0{\bm e}\sb{1,c}
-R\sb 1{\bm e}\sb{2,c}
-R\sb 2{\bm e}\sb{3,c}
-\zeta(\zeta+R\sb 1)\p\sb{c}{\bm e}\sb{3,c}
+R\sb 0\p\sb{x}(\zeta{\bm e}\sb{3,c}+\bm\upsilon)
+J{\bm N},
\end{equation}
where
$c=c(t)=c\sb\star+\eta(t)$,
$\zeta=\zeta(t)$,
and the nonlinear terms $R\sb j(t)$ are given by (\ref{def-r}).
We set
\[
\bm\omega(x,t)=e\sp{\mu x}\bm\upsilon(x,t),
\qquad
{\bm e}\sp{\mu}\sb{j,c}(x)=e\sp{\mu x}{\bm e}\sb{j,c}(x),
\qquad
c\in[c\sb\star,c\sb\star+\eta\sb 1],
\quad
j=1,\,2,\,3,
\]
and consider $A\sp{\mu}\sb{c}$
given by (\ref{conjugate}).
Equation (\ref{evolution-upsilon}) takes the following form:
\begin{equation}\label{proof-problem1}
\dot{\bm\omega}-A\sp{\mu}\sb{c}\bm\omega={\bm G},
\end{equation}
where
\begin{equation}
{\bm G}(x,t)=-R\sb 0{\bm e}\sp{\mu}\sb{1,c}
-R\sb 1{\bm e}\sp{\mu}\sb{2,c}-R\sb 2{\bm e}\sp{\mu}\sb{3,c}
-\zeta(\zeta+R\sb 1)\p\sb{c}{\bm e}\sp{\mu}\sb{3,c}
+R\sb 0(\p\sb x-\mu)(\zeta{\bm e}\sp{\mu}\sb{3,c}+\bm\omega)
+e\sp{\mu x}J{\bm N}.
\end{equation}
As follows from (\ref{proof-problem1}),
\[
\p\sb{t}({Q\sp{\mu}\sb{c\sb\star}\bm\omega})
=Q\sp{\mu}\sb{c\sb\star}\dot{\bm\omega}
=Q\sp{\mu}\sb{c\sb\star}
\left(
A\sp{\mu}\sb{c}\bm\omega+{\bm G}
\right)
=A\sp{\mu}\sb{c\sb\star}Q\sp{\mu}\sb{c\sb\star}\bm\omega
+Q\sp{\mu}\sb{c\sb\star}(A\sp{\mu}\sb{c}-A\sp{\mu}\sb{c\sb\star})\bm\omega
+Q\sp{\mu}\sb{c\sb\star}{\bm G}.
\]
We may write $Q\sp{\mu}\sb{c\sb\star}\bm\omega$ as follows:
\begin{equation}\label{345}
Q\sp{\mu}\sb{c\sb\star}\bm\omega(t)
=
\int\sb 0^t e\sp{A\sp{\mu}\sb{c\sb\star}(t-s)}\mathfrak{G}(s)\,ds,
\end{equation}
where
\begin{equation}\label{def-g}
\mathfrak{G}(x,t)=
Q\sp{\mu}\sb{c\sb\star}(A\sp{\mu}\sb{c}-A\sp{\mu}\sb{c\sb\star})\bm\omega(x,t)
+Q\sp{\mu}\sb{c\sb\star}{\bm G}(x,t).
\end{equation}
Using the dissipative estimate given by (\ref{smoothing}),
we have:
\begin{eqnarray}\label{help-bounds}
\norm{Q\sp{\mu}\sb{c\sb\star}\bm\omega(t)}\sb{H\sp{1}}
\le C
\int\sb 0\sp{t}
(t-s)\sp{-1/2}
e\sp{-b(t-s)}
\norm{\mathfrak{G}(s)}\sb{L\sp{2}}
\,ds
\\
\le C
e\sp{-b t/2}\sup\sb{s\in[0,t]}
e^{b s/2}
\norm{\mathfrak{G}(s)}\sb{L\sp{2}}
\int\sb 0\sp{t}
(t-s)\sp{-1/2}
e\sp{-b(t-s)/2}
\,ds
\\
\le C
\sup\sb{s\in[0,t]}
e^{b s/2}
\norm{\mathfrak{G}(s)}\sb{L\sp{2}}.
\end{eqnarray}
Since
$\bm\omega
=Q\sp{\mu}\sb{c}\bm\omega
=Q\sp{\mu}\sb{c\sb\star}\bm\omega
+(Q\sp{\mu}\sb{c}-Q\sp{\mu}\sb{c\sb\star})\bm\omega$,
we have
\[
\norm{\bm\omega}\sb{H\sp{1}}
\le\norm{
Q\sp{\mu}\sb{c\sb\star}\bm\omega}\sb{H\sp{1}}
+\abs{\eta}
\sup\sb{c\in[c\sb\star,c\sb\star+\eta\sb 1]}
\norm{\p\sb{c}Q\sp{\mu}\sb{c}}\sb{H\sp 1\to H\sp 1}
\norm{\bm\omega}\sb{H\sp{1}}
\le\norm{
Q\sp{\mu}\sb{c\sb\star}\bm\omega}\sb{H\sp{1}}
+
\frac 1 2
\norm{\bm\omega}\sb{H\sp{1}},
\]
where we used the inequality (\ref{eta-small-q}).
It follows that 
$\norm{\bm\omega}\sb{H\sp{1}}
\le 2\norm{Q\sp{\mu}\sb{c\sb\star}\bm\omega}\sb{H\sp{1}}.
$
Hence, we have:
\begin{eqnarray}\label{help-bounds-05}
\norm{\bm\omega(t)}\sb{H\sp{1}}
\le C
e\sp{-b t/2}
\sup\sb{s\in[0,t]}e\sp{b s/2}
\norm{\mathfrak{G}(s)}\sb{L\sp{2}}.
\end{eqnarray}
We now need the bound on
$\norm{\mathfrak{G}}\sb{L\sp{2}}$:
\begin{equation}\label{GG}
\norm{\mathfrak{G}}\sb{L\sp{2}}
\le
\norm{Q\sp{\mu}\sb{c\sb\star}(A\sp{\mu}\sb{c}-A\sp{\mu}\sb{c\sb\star})\bm\omega}
\sb{L\sp{2}}
+\norm{Q\sp{\mu}\sb{c\sb\star}\bm G}\sb{L\sp{2}}.
\end{equation}
We estimate the first term in the right-hand side of (\ref{GG})
as follows:
\begin{equation}\label{help-bounds1}
\norm{
Q\sp{\mu}\sb{c\sb\star}(A\sp{\mu}\sb{c}-A\sp{\mu}\sb{c\sb\star})\bm\omega(t)
}\sb{L\sp{2}}
\le
\norm{
Q\sp{\mu}\sb{c\sb\star}(A\sp{\mu}\sb{c}-A\sp{\mu}\sb{c\sb\star})
}\sb{H\sp{1}\to L\sp{2}}\norm{\bm\omega(t)}\sb{H\sp{1}}
\le C\abs{\eta}\norm{\bm\omega(t)}\sb{H\sp{1}}.
\end{equation}
Since ${\bm e}\sp{\mu}\sb{j,c}$,
$1\le j\le 3$, depend continuously on $c$
while
$Q\sp{\mu}\sb{c\sb\star}{\bm e}\sp{\mu}\sb{j,c\sb\star}=0$,
there are bounds
$
\norm{Q\sp{\mu}\sb{c\sb\star}{\bm e}\sp{\mu}\sb{j,c}}\sb{H\sp{1}}
\le C\abs{\eta}$.
This allows to derive the following bound for the second term
in the right-hand side of (\ref{GG}):
\begin{eqnarray*}
\norm{Q\sp{\mu}\sb{c\sb\star}{\bm G}}\sb{L\sp{2}}
\le C
\left(
\abs{\eta}
\sup\sb{0\le j\le 2}
\abs{R\sb{j}}
+\abs{\zeta}\abs{\zeta+R\sb 1}
+\abs{R\sb 0}
(\abs{\zeta}+\norm{\bm\omega}\sb{H\sp{1}})
+\norm{J{\bm N}}\sb{L\sb\mu\sp{2}}
\right)
\\
\le
C\left(
\zeta^2
+
(\abs{\eta}+\abs{\zeta}+\norm{\bm\omega}\sb{H\sp{1}})
\sup\sb{0\le j\le 2}\abs{R\sb j}
+\norm{J{\bm N}}\sb{L\sb\mu\sp{2}}
\right).
\end{eqnarray*}
Using the representation (\ref{def-r})
and the inclusions
$\p\sb{c}{\bm e}\sb{3,c}\in H\sb{\mu}\sp\infty(\R)$,
${\bm g}\sb i\in H\sb{-\mu}\sp\infty(\R)$,
$\p\sb{c}{\bm g}\sb i\in H\sb{-\mu}\sp\infty(\R)$,
we obtain the following estimates on $R\sb j$:
\begin{equation}
\abs{R\sb j(\eta,\zeta,\bm\upsilon)}
\le
C\left(
\zeta^2
+\abs{\zeta}\norm{\bm\upsilon}\sb{H\sb\mu\sp 1}
+\norm{J{\bm N}}\sb{L\sb\mu\sp 2}
\right),
\qquad
j=0,\,1,\,2.
\label{bound-on-r}
\end{equation}
Taking into account (\ref{bound-on-r}), we get:
\begin{eqnarray}
\norm{Q\sp{\mu}\sb{c\sb\star}\bm G}\sb{L\sb\mu\sp{2}}
\le C\left(
\zeta^2
+(\abs{\eta}+\abs{\zeta}+\norm{\bm\omega}\sb{H\sp{1}})
(\zeta^2+\abs{\zeta}\norm{\bm\omega}\sb{H\sp{1}}
+\norm{J{\bm N}}\sb{L\sb\mu\sp{2}})
+\norm{J{\bm N}}\sb{L\sb\mu\sp{2}}
\right)
\nonumber
\\
\le C\left(
\zeta^2
+(\abs{\eta}+\norm{\bm\omega}\sb{H\sp{1}})
\abs{\zeta}\norm{\bm\omega}\sb{H\sp{1}}
+\norm{J{\bm N}}\sb{L\sb\mu\sp{2}}
\right).
\label{help-bounds2}
\end{eqnarray}
In the last inequality, we used the uniform boundedness
of $\abs{\eta}$, $\abs{\zeta}$,
and $\norm{\bm\omega}\sb{H\sp{1}}$
that follows from (\ref{uniform-bounds-assumed}).

Summing up (\ref{help-bounds1}) and (\ref{help-bounds2}),
we obtain the following bound on $\norm{\mathfrak{G}}\sb{L\sb\mu\sp{2}}$:
\begin{equation}\label{bound-g-l2}
\norm{\mathfrak{G}}\sb{L\sb\mu\sp{2}}
\le C\left[
\zeta^2
+(\abs{\eta}+\abs{\zeta})\norm{\bm\omega}\sb{H\sp{1}}
+\norm{J{\bm N}}\sb{L\sb\mu\sp{2}}
\right].
\end{equation}
Using the integral representation
for the nonlinearity (\ref{def-n}),
\begin{equation}\label{taylor-jn}
J{\bm N}
=\p\sb x[f(\bm\phi\sb{c}+\bm\rho)-f(\bm\phi\sb{c})-f'(\bm\phi\sb{c})\bm\rho]
=\p\sb x
\left[
\frac{\bm\rho^2}{2}\int\sb 0\sp 1
(1-s)^2 f''(\bm\phi\sb{c}+s\bm\rho)\,ds
\right],
\end{equation}
we obtain the bound
\[
\norm{J{\bm N}}\sb{L\sb\mu\sp{2}}
\le C
\norm{\bm\rho}\sb{H\sb\mu\sp{1}}\norm{\bm\rho}\sb{H\sp{1}}
\le C
\left(
\abs{\zeta}
\norm{{\bm e}\sb{3,c}}\sb{H\sb\mu\sp 1}
+\norm{\bm\upsilon}\sb{H\sb\mu\sp{1}}
\right)
\norm{\bm\rho}\sb{H\sp{1}},
\]
with the constant $C$
that depends on $\norm{\bm\phi\sb c}\sb{H\sp 1}$
and on the bounds on $f''(z)$ and $f'''(z)$
for
$\abs{z}\le\norm{{\bm u}}\sb{L\sp\infty}$,
which is bounded by
$2\norm{\bm\phi\sb{c\sb\star}}\sb{H\sp 1}$.
This bound allows to rewrite (\ref{bound-g-l2}) as
\begin{eqnarray}
\norm{\mathfrak{G}}\sb{L\sb\mu\sp{2}}
\le C\left[
\zeta^2
+(\abs{\eta}+\abs{\zeta}+\norm{\bm\rho}\sb{H\sp{1}})\norm{\bm\omega}\sb{H\sp{1}}
+\abs{\zeta}\norm{\bm\rho}\sb{H\sp{1}}
\right]
\le C\left[
g\sb 0
+g\sb 1
\norm{\bm\omega}\sb{H\sp{1}}
\right],
\end{eqnarray}
where
\begin{equation}\label{def-g-0-1}
g\sb 0(t)
=\zeta^2(t)+\abs{\zeta(t)}\norm{\bm\rho(t)}\sb{H\sp{1}},
\qquad
g\sb 1(t)
=\abs{\eta(t)}+\abs{\zeta(t)}+\norm{\bm\rho(t)}\sb{H\sp{1}}.
\end{equation}
Thus,
(\ref{help-bounds-05})
could be written as
\begin{eqnarray}\label{help-bounds-06}
2 e\sp{b t/2}
\norm{\bm\omega(t)}\sb{H\sp{1}}
\le C\sb 2
\sup\sb{s\in[0,t]}e\sp{b s/2}
\left[
g\sb 0(s)
+g\sb 1(s)\norm{\bm\omega(s)}\sb{H\sp{1}}
\right],
\end{eqnarray}
for some $C\sb 2>0$.
Since the right-hand side is monotonically increasing with $t$,
we also have
\begin{eqnarray}\label{help-bounds-08}
\sup\sb{s\in[0,t]}
2e\sp{b s/2}
\norm{\bm\omega(s)}\sb{H\sp{1}}
\le
C\sb 2
\sup\sb{s\in[0,t]}
e\sp{b s/2}
\left[
g\sb 0(s)
+
g\sb 1(s)
\norm{\bm\omega(s)}\sb{H\sp{1}}
\right].
\end{eqnarray}
The function $g\sb 1$ from (\ref{def-g-0-1})
satisfies
$
C\sb 2\sup\sb{s\in[0,t]}g\sb 1(s)<1
$
(this follows from
the assumptions (\ref{eta-zeta-delta-small})
and
(\ref{uniform-bounds-assumed})),
and therefore
\[
\norm{\bm\omega(t)}\sb{H\sp{1}}
\le C\sb 2
e\sp{-b t/2}
\sup\sb{s\in[0,t]}e\sp{b s/2}g\sb 0(s)
\le C\sb 2
\sup\sb{s\in[0,t]}
\left[
\zeta^2(s)+\abs{\zeta(s)}\norm{\bm\rho(s)}\sb{H\sp{1}}
\right].
\]
Since $\bm\omega=e\sp{\mu x}\bm\upsilon$,
the last inequality yields (\ref{bound-v-mu0}).
\end{proof}

\section{Nonlinear estimates}
\label{sect-nonlinear-estimates}

Now we close the estimates using the bounds on
$\norm{\bm\rho}\sb{H\sp 1}$ (Lemma~\ref{lemma-energy}) and on
$\norm{\bm\upsilon}\sb{H\sb\mu\sp{1}}$ (Lemma
\ref{lemma-dispersion}) from the previous section.

We assume that
$\eta\sb 1>0$, $\zeta\sb 1>0$,
and $\delta\sb 1>0$
are sufficiently small:
not larger than in
Proposition~\ref{prop-well},
satisfy the bounds
(\ref{delta-small-1}), (\ref{eta-small-q}), and (\ref{eta-zeta-delta-small}),
and also that $\zeta\sb 1$
satisfies
\begin{equation}\label{zeta-small-1}
\zeta\sb 1
<\frac{1}{3 \max(1,C\sb 1)C\sb 2}.
\end{equation}
Define
\begin{equation}\label{c3-c4}
C\sb 3=2 C\sb 1,
\qquad
C\sb 4=2 C\sb 2\max(1,C\sb 3),
\end{equation}
with $C\sb 1$ and $C\sb 2$ as in Lemmas~\ref{lemma-energy} and
\ref{lemma-dispersion}.
Choosing smaller values of $\eta\sb 1$ and $\zeta\sb 1$ if necessary,
we may assume that
\begin{eqnarray}
&&
C\sb 3\left(
\zeta\sb 1
+2\eta\sb 1
+(\mathscr{N}\sb{c\sb\star+\eta\sb 1}-\mathscr{N}\sb{c\sb\star})\sp{1/2}
\right)
<\delta\sb 1,
\label{cc-1}
\\
&&
C\sb 4\left(
\zeta\sb 1^2
+2\eta\sb 1\zeta\sb 1
+
\zeta\sb 1
(\mathscr{N}\sb{c\sb\star+\eta\sb 1}
-\mathscr{N}\sb{c\sb\star})\sp{1/2}
\right)
<\delta\sb 1.
\label{cc-2}
\end{eqnarray}

Define
\begin{eqnarray}
&&
\eta\sb M(t)=\sup\sb{0\le s\le t}\eta(s),
\label{def-eta-m}
\\
&&
\zeta\sb M(t)=\sup\sb{0\le s\le t}\abs{\zeta(s)}.
\label{def-zeta-m}
\end{eqnarray}

\begin{proposition}\label{prop-5}
Assume that the initial data $\eta\sb 0>0$ and $\zeta\sb 0$
are such that the following inequalities are satisfied:
\begin{equation}\label{init}
\eta\sb 0\in(0,\eta\sb 1),
\qquad
\abs{\zeta\sb 0}<\zeta\sb 1,
\qquad
\norm{\bm\rho\sb 0}\sb{H\sp{1}}
<\min(\eta\sb 1,\delta\sb 1).
\end{equation}
Then for $0\le t<T\sb 1$
the functions $\bm\rho(t)$, $\bm\upsilon(t)$ satisfy the bounds
\begin{eqnarray}
\label{bound-rho}
&&
\norm{\bm\rho(t)}\sb{H\sp{1}}
\le C\sb 3
\left[
\zeta\sb M(t)+\mathscr{Y}(\eta\sb M(t))
\right],
\\
\label{bound-v-mu}
&&\norm{\bm\upsilon(t)}\sb{H\sb\mu\sp{1}}
\le C\sb 4
\left[
\zeta\sb M(t)^2+\zeta\sb M(t)\mathscr{Y}(\eta\sb M(t))
\right],
\end{eqnarray}
where $C\sb 3$, $C\sb 4$ are defined by (\ref{c3-c4}),
$
\mathscr{Y}(\eta)
=
\norm{\bm\rho\sb 0}\sb{H\sp{1}}
+\norm{\bm\rho\sb 0}\sb{H\sp{1}}\sp{1/2}\abs{\eta-\eta\sb 0}^{1/2}
+\abs{\mathscr{N}\sb{c\sb\star+\eta}
-\mathscr{N}\sb{c\sb\star+\eta\sb 0}}\sp{1/2}
$
is introduced in (\ref{def-y}),
and $\eta\sb{M}$, $\zeta\sb{M}$
are defined in {\rm (\ref{def-eta-m})} and {\rm (\ref{def-zeta-m})}.
\end{proposition}

\begin{proof}
Let
\[
S=\{t\in[0,T\sb 1)\sothat
\norm{\bm\rho(t)}\sb{H\sp{1}}<\delta\sb 1
\}.
\]
$S$ is nonempty
since
$\norm{\bm\rho(0)}\sb{H\sp{1}}<\delta\sb 1$
by (\ref{init}).
According to
Proposition~\ref{prop-well}
and representation (\ref{def-rho}),
$\norm{\bm\rho(t)}\sb{H\sp{1}}$
is a continuous function of $t$.
Since the inequality in the definition of $S$
is sharp, $S$ is an open subset of $[0,T\sb 1)$.
Let us assume that
$T\sb 2\in(0,T\sb 1)$ is such that
\begin{equation}\label{assumption}
\norm{\bm\rho(t)}\sb{H\sp{1}}
<\delta\sb 1,
\qquad
0\le t<T\sb 2.
\end{equation}
It is enough to prove that $T\sb 2\in S$
(then the connected subset of $S$ that contains $t=0$
is both open and closed in $[0,T\sb 1)$
and hence coincides with $[0,T\sb 1)$).
Since $\norm{\bm\upsilon(t)}\sb{H\sb\mu\sp{1}}<\delta\sb 1$
for $0\le t<T\sb 1$,
both Lemma~\ref{lemma-energy} and Lemma~\ref{lemma-dispersion}
are applicable
for $t\le T\sb 2$.
The estimate (\ref{bound-rho0}) on $\norm{\bm\rho(t)}\sb{H\sp{1}}$
together with the estimate (\ref{bound-v-mu0})
on $\norm{\bm\upsilon(t)}\sb{H\sb\mu\sp{1}}$
give
\[
\norm{\bm\rho(t)}\sb{H\sp{1}}
\le C\sb 1\left(
\mathscr{Y}(\eta(t))+\abs{\zeta(t)}
+\norm{\bm\upsilon(t)}\sb{H\sb\mu\sp{1}}\right)
\le C\sb 1\Big(
\mathscr{Y}(t)+\abs{\zeta(t)}
+C\sb 2
\sup\sb{s\in[0,t]}
\left[\zeta^2+\abs{\zeta}\norm{\bm\rho}\sb{H\sp{1}}\right]
\Big).
\]
For $0\le t\le T\sb 2$,
define
$M(t)=\sup\sb{s\in[0,t]}\norm{\bm\rho(s)}\sb{H\sp{1}}$.
We have:
\[
M(t)
\le C\sb 1
\Big(
\sup\sb{s\in[0,t]}\left(\mathscr{Y}(\eta(s))+\abs{\zeta(s)}\right)
+
C\sb 2
\sup\sb{s\in[0,t]}
\left[
\zeta^2(s)+\abs{\zeta(s)}M(t)
\right]
\Big).
\]
We carry the term $C\sb 1C\sb 2\abs{\zeta}M(t)$
to the left-hand side of the inequality,
taking into account that
$C\sb 1 C\sb 2\abs{\zeta(t)}
\le C\sb 1C\sb 2\zeta\sb 1\le\frac{1}{3}$
for all $0\le t<T\sb 1$
by (\ref{zeta-small-1}).
This results in the following relation:
\[
\norm{\bm\rho(t)}\sb{H\sp{1}}
\le M(t)\le \frac{3}{2}C\sb 1
\Big(
\sup\sb{s\in[0,t]}
\left(\mathscr{Y}(\eta(s))+\abs{\zeta(s)}\right)
+C\sb 2
\sup\sb{s\in[0,t]}\zeta^2(s)
\Big).
\]
Since
$C\sb 2\zeta^2\le C\sb 2\zeta\sb 1\abs{\zeta}
\le\abs{\zeta}/3$
by (\ref{zeta-small-1}),
we obtain:
\[
\norm{\bm\rho(t)}\sb{H\sp{1}}
\le \frac 3 2
C\sb 1\sup\sb{s\in[0,t]}
\Big(\mathscr{Y}(\eta(s))+\frac 4 3 \abs{\zeta(s)}\Big)
\le C\sb 3\sup\sb{s\in[0,t]}\left(\mathscr{Y}(\eta(s))+\abs{\zeta(s)}\right),
\qquad
t\in[0,T\sb 2],
\]
with $C\sb 3=2 C\sb 1$.
This proves (\ref{bound-rho}) for $t\in[0,T\sb 2]$.
It then follows that
\[
\norm{\bm\rho(T\sb 2)}\sb{H\sp{1}}
\le
C\sb 3
\left[
\zeta\sb 1
+\mathscr{Y}(\eta\sb 1)
\right]
\le
C\sb 3
\left[
\zeta\sb 1
+2\eta\sb 1
+\left(
\mathscr{N}\sb{c\sb\star+\eta\sb 1}
-\mathscr{N}\sb{c\sb\star}
\right)\sp{1/2}
\right]
<\delta\sb 1,
\]
where we took into account
the definition of $\mathscr{Y}(\eta)$
in (\ref{def-y}),
the bound
$\norm{\bm\rho\sb 0}\sb{H\sp 1}<\eta\sb 1$
from (\ref{init}),
and the inequality
(\ref{cc-1}).
Hence, $T\sb 2\in S$.
It follows that $S$ coincides with $[0,T\sb 1)$.


Using  the bound (\ref{bound-rho})
in (\ref{bound-v-mu0})
and recalling the definition of $C\sb 4$ in (\ref{c3-c4}),
we derive the bound (\ref{bound-v-mu})
on $\norm{\bm\upsilon(t)}\sb{H\sb\mu\sp{1}}$.
\end{proof}

\begin{corollary}\label{corollary-mu-2}
Assume that conditions of Proposition~\ref{prop-5}
are satisfied.
If $\eta\sb 1>0$ and $\zeta\sb 1>0$
were chosen sufficiently small,
then there exists a constant $C\sb 5>0$ 
so that for $0\le t<T\sb 1$
the function $\bm\upsilon(t)$
satisfies the bound
\begin{equation}\label{bound-v-mu-2}
\norm{\bm\upsilon(t)}\sb{H\sb{\mu/2}\sp{1}}
\le C\sb 5
\left[
\zeta\sb M^2(t)
+\zeta\sb M(t)
\mathscr{Y}(\eta\sb M(t))
\right],
\end{equation}
where $\eta\sb{M}$, $\zeta\sb{M}$
are defined in {\rm (\ref{def-eta-m})}, {\rm (\ref{def-zeta-m})}.
\end{corollary}

\begin{proof}
The bound (\ref{bound-v-mu-2})
is proved in the same way as (\ref{bound-v-mu}).
We may need to take smaller values of
$\eta\sb 1$ and $\zeta\sb 1$
so that
Lemmas~\ref{lemma-energy} and~\ref{lemma-dispersion}
become applicable for the new exponential weight.
Note that the exponential weight
does not enter the definition (\ref{def-y})
of the function $\mathscr{Y}(\eta)$.
\end{proof}

\begin{lemma}\label{lemma-bounds-on-r}
Assume that the bounds {\rm (\ref{bound-v-mu})}
and {\rm (\ref{bound-v-mu-2})}
are satisfied for $0\le t<T\sb 1$.
Then there exists $C\sb 6>0$
so that
the terms $R\sb 1$ and $R\sb 2$
defined in {\rm (\ref{def-r})}
satisfy
for $0\le t<T\sb 1$
the bounds
\begin{equation}
\label{bound-r}
\abs{R\sb{j}(\eta,\zeta,\bm\upsilon)}
\le
C\sb 6\zeta\sb M^2,
\qquad
j=1,\,2.
\end{equation}
\end{lemma}

\begin{proof}
By (\ref{bound-on-r}),
\begin{equation}\label{rj}
\abs{R\sb{j}(\eta,\zeta,\bm\upsilon)}
\le
C\left(
\zeta^2+\abs{\zeta}\norm{\bm\upsilon}\sb{H\sb\mu\sp{1}}
+\norm{J{\bm N}}\sb{L\sb\mu\sp{2}}
\right),
\qquad
j=1,\,2.
\end{equation}
According to (\ref{bound-v-mu}),
the second term in the right-hand side of (\ref{rj})
is bounded by
$C\zeta^2$
as long as $\eta\in(0,\eta\sb 1)$
and
$\abs{\zeta}\le\zeta\sb 1$.
We now need a bound on $\norm{J{\bm N}}\sb{L\sb\mu\sp{2}}$.
Using the representation (\ref{taylor-jn})
for the nonlinearity,
we obtain the bounds
\begin{equation}\label{help-1}
\norm{J{\bm N}}\sb{L\sb\mu\sp{2}}
\le C
\norm{\bm\rho}\sb{H\sb{\mu/2}\sp{1}}^2
\le C
\left(
\zeta^2\norm{{\bm e}\sb{3,c}}\sb{H\sb{\mu/2}\sp{1}}^2
+\norm{\bm\upsilon}\sb{H\sb{\mu/2}\sp{1}}^2
\right).
\end{equation}
The constant depends on
$\norm{\bm\phi\sb c}\sb{H\sp 1}$
and on the bounds on $f''(z)$ and $f'''(z)$
for
$\abs{z}\le\norm{{\bm u}}\sb{L\sp\infty}$,
which is bounded by
$2\norm{\bm\phi\sb{c\sb\star}}\sb{H\sp 1}$.
As follows from (\ref{bound-v-mu-2}),
\begin{equation}\label{help-15}
\norm{\bm\upsilon(t)}\sb{H\sb{\mu/2}\sp{1}}
\le
C\sb 5(\zeta\sb 1+\mathscr{Y}(\eta\sb 1))
\zeta\sb M(t).
\end{equation}
Using this bound in (\ref{help-1}), we get
$\norm{J{\bm N}}\sb{L\sb\mu\sp{2}}\le C\zeta\sb M^2$.
The bound (\ref{bound-r}) follows.
\end{proof}

\section{Choosing the initial perturbation}
\label{sect-general-case}

In this section,
we show how to choose the initial perturbation
that indeed leads to the instability
and conclude the proof of Theorem~\ref{main-theorem}.

We choose
$\eta\sb 1>0$, $\zeta\sb 1>0$, and $\delta\sb 1>0$
small enough
so that the inequalities
(\ref{delta-small-1}), (\ref{eta-small-q}), (\ref{eta-zeta-delta-small}),
are satisfied,
and so that
Lemmas~\ref{lemma-energy} and~\ref{lemma-dispersion}
apply to both exponential weights $\mu$ and $\mu/2$.
Taking $\eta\sb 1>0$, $\zeta\sb 1>0$
smaller if necessary,
we may assume that the conditions
(\ref{zeta-small-1}),
(\ref{cc-1}), and (\ref{cc-2})
are satisfied,
and moreover that
\begin{equation}\label{zeta1-small}
C\sb 6\zeta\sb 1<1/2,
\end{equation}
where $C\sb 6>0$ is from Lemma~\ref{lemma-bounds-on-r}.

Let
\begin{equation}\label{def-anti-lambda}
\lambda(\eta)=\lambda\sb{c\sb\star+\eta},
\qquad
\Lambda(\eta)=\int\sb 0\sp{\eta}\lambda(\eta')\,d\eta'.
\end{equation}
Let us recall that, according to (\ref{def-eta1}),
we assume that there exists $\eta\sb 1>0$
so that $\mathscr{N}\sb{c}'<0$ and is nonincreasing
for $c\sb\star<c\le c\sb\star+\eta\sb 1$.
Thus,
we assume that $\lambda(\eta)>0$ for $0<\eta\le\eta\sb 1$
(according to (\ref{def-lambda}),
$\mathscr{N}\sb{c}'$ and $\lambda\sb{c}$
are of opposite sign).

\begin{lemma}\label{lemma-3e}
One can choose $\eta\sb 1>0$ sufficiently small
so that
for $0<\eta\le\eta\sb 1$
one has
\begin{equation}\label{3e}
3 C\sb 6 e^{2C\sb 6\eta}\Lambda(\eta)
<\lambda(\eta).
\end{equation}
\end{lemma}

\begin{proof}
By (\ref{def-lambda}),
$
\lambda\sb c
=-\frac{\mathscr{N}\sb c'}{B\sb c},
$
where
\begin{equation}
B\sb c=\langle\bm\phi\sb c,{\bm e}\sb{3,c}\rangle.
\end{equation}
Since $B\sb{c\sb\star}>0$ by (\ref{fourth}),
we may assume that $\eta\sb 1>0$
is small enough so that
\begin{equation}\label{b2b2}
B\sb{c\sb\star}/2\le B\sb c\le 2B\sb{c\sb\star},
\qquad
c\in[c\sb\star,c\sb\star+\eta\sb 1].
\end{equation}
According to Theorem~\ref{main-theorem},
$\mathscr{N}\sb{c}'<0$ and is nonincreasing
for $c\in(c\sb\star,c\sb\star+\eta\sb 1)$.
Therefore, using inequalities (\ref{b2b2}), we obtain:
\[
\Lambda(\eta)
=\int\limits\sb{c\sb\star}\sp{c\sb\star+\eta}\lambda\sb{c}\,d c
=\int\limits\sb{c\sb\star}\sp{c\sb\star+\eta}
\frac{-\mathscr{N}\sb c'}{B\sb c}\,d c
\le
-\frac{2\eta\mathscr{N}\sb{c\sb\star+\eta}'}{B\sb{c\sb\star}}
\le
4\eta\lambda(\eta),
\qquad
0\le\eta\le\eta\sb 1,
\]
where $\lambda(\eta)>0$ for $0<\eta\le\eta\sb 1$.
We take $\eta\sb 1>0$ so small that
$
12\eta\sb 1 C\sb 6 e^{2 C\sb 6\eta\sb 1}<1;
$
then (\ref{3e}) is satisfied.
\end{proof}

Taking $\eta\sb 1>0$ smaller if necessary,
we may assume that
Lemma~\ref{lemma-nd} is satisfied and that
\begin{equation}\label{eta1-small}
\lambda(\eta)/C\sb 6<\zeta\sb 1.
\end{equation}

\begin{remark}
The inequality (\ref{eta1-small})
ensures that $\eta(t)$ reaches $\eta\sb 1$
prior to $\zeta(t)$ reaching $\zeta\sb 1$
(see Lemma~\ref{lemma-zeta-le-Lambda-2} and Figure~\ref{fig-eta-zeta-mono}).
\end{remark}

Since $\Lambda(\eta)=\os(\eta)$,
we may also assume that $\eta\sb 1>0$ is small enough so that
\begin{equation}\label{c7-lambda}
K(\eta\sb 1,\zeta\sb 1)\Lambda(\eta\sb 1)\le\kappa\eta\sb 1/2,
\end{equation}
where the function
$K(\eta\sb 1,\zeta\sb 1)$ is defined below in (\ref{def-c7})
and
$\kappa>0$ is from Lemma~\ref{lemma-nd}.

\begin{lemma}\label{lemma-init-data}
For any
$\delta\in(0,\min(\eta\sb 1,\delta\sb 1))$,
one can choose the initial data
$
\eta\sb 0\in(0,\eta\sb 1),
$
$
\zeta\sb 0\in (0,\zeta\sb 1)
$
so that the following estimates are satisfied:
\begin{equation}\label{c-1}
\norm{
\zeta\sb 0{\bm e}\sb{3,c\sb\star+\eta\sb 0}
}\sb{H\sp 1}
<\min(\eta\sb 1,\delta\sb 1),
\end{equation}
\begin{equation}\label{c-2}
\norm{
(\bm\phi\sb{c\sb\star+\eta\sb 0}
+\zeta\sb 0{\bm e}\sb{3,c\sb\star+\eta\sb 0})
-\bm\phi\sb{c\sb\star}}
\sb{H\sp{1}\cap H\sb\mu\sp 1}
<\delta<\min(\eta\sb 1,\delta\sb 1),
\end{equation}
\begin{equation}\label{c-3}
\zeta\sb 0<\Lambda(\eta\sb 0).
\end{equation}
\end{lemma}

\begin{proof}
Pick $\eta\sb 0\in(0,\eta\sb 1)$
so that
\begin{equation}\label{eta0-small-enough}
\norm{\bm\phi\sb{c\sb\star+\eta\sb 0}-\bm\phi\sb{c\sb\star}}
\sb{H\sp 1\cap H\sb\mu\sp 1}
<\delta/2.
\end{equation}
For given $\eta\sb 0>0$,
we take $\zeta\sb 0\in(0,\zeta\sb 1)$
small enough so that
\begin{equation}\label{zeta0-small-enough}
\zeta\sb 0
\norm{{\bm e}\sb{3,c\sb\star+\eta\sb 0}}\sb{H\sp 1\cap H\sb\mu\sp 1}
<\delta/2.
\end{equation}
Note that $\norm{{\bm e}\sb{3,c\sb\star+\eta\sb 0}}\sb{H\sp 1}$
for $\eta\sb 0>0$
is finite
by Lemma~\ref{lemma-is}.
Inequality
(\ref{zeta0-small-enough}) implies that
(\ref{c-1}) is satisfied.
Together with
(\ref{eta0-small-enough}),
it also guarantees that (\ref{c-2}) holds.
We then require that $\zeta\sb 0>0$
be small enough so that the inequality
(\ref{c-3}) takes place.
\end{proof}

We rewrite the two last equations from the system (\ref{xi-eta-zeta-system}):
\begin{equation}\label{system-degenerate}
\left\{
\begin{array}{l}
\dot\eta=\zeta+R\sb 1(\eta,\zeta,\bm\upsilon),
\\~\\
\dot\zeta=\lambda(\eta)\zeta+R\sb 2(\eta,\zeta,\bm\upsilon).
\end{array}
\right.
\end{equation}

\begin{lemma}\label{lemma-zeta-le-Lambda-2}
For $0\le t<T\sb 1$,
with $T\sb 1>0$ as in Proposition~\ref{prop-well},
\begin{eqnarray}
&&
\dot\eta\ge\zeta\sb 0/2,
\qquad
\dot\zeta\ge 0,
\label{monotone}
\\
&&
\zeta\sb 0\le\zeta(t)< 3 e^{2C\sb 6\eta(t)}\Lambda(\eta(t)).
\label{zeta-bounds}
\end{eqnarray}
\end{lemma}

\begin{proof}
According to Proposition~\ref{prop-well},
the trajectory $(\eta(t),\zeta(t))$
that starts at $(\eta\sb 0,\zeta\sb 0)$
satisfies the inequalities
$\eta(t)<\eta\sb 1$ and $\zeta(t)<\zeta\sb 1$
for $0\le t<T\sb 1$.
We define the region $\Omega\subset\R\sb{+}\times\R\sb{+}$
by
\begin{equation}\label{def-omega}
\Omega=\{
(\eta,\zeta)\sothat
\zeta\sb 0\le\zeta\le\lambda(\eta)/C\sb 6,
\ \eta\sb 0\le\eta\le\eta\sb 1
\}.
\end{equation}
Define $T\sb{\Omega}\in\R\sb{+}\cup\{+\infty\}$ by
\begin{equation}\label{def-t22}
T\sb{\Omega}=\sup\{\,t\in[0,T\sb 1)\,\sothat
\ (\eta(t),\zeta(t))\in\Omega,
\quad \dot\zeta(t)\ge 0
\,\}.
\end{equation}
Let us argue that
$T\sb{\Omega}>0$.
At $t=0$,
$(\eta(0),\zeta(0))=(\eta\sb 0,\zeta\sb 0)\in\Omega$.
From (\ref{system-degenerate}),
we compute:
$\dot\eta(0)\ge\zeta\sb 0-C\sb 6\zeta\sb 0^2>0$,
where we applied the bounds (\ref{bound-r})
and the inequality $C\sb 6\zeta\sb 0<1/2$ 
that follows from (\ref{zeta1-small})
and the choice $\zeta\sb 0<\zeta\sb 1$.
Similarly,
$\dot\zeta(0)\ge\lambda(\eta\sb 0)\zeta\sb 0-C\sb 6\zeta\sb 0^2>0$
due to the inequality $C\sb 6\zeta\sb 0<\lambda(\eta\sb 0)$
that follows from
(\ref{c-3}) and (\ref{3e}).
Therefore,
$(\eta(t),\zeta(t))\in\Omega$
and $\dot\zeta(t)>0$
for times $t>0$ from a certain open neighborhood of $t=0$,
proving that $T\sb{\Omega}>0$.

The monotonicity of $\zeta(t)$
for $t<T\sb{\Omega}$ implies that
$\zeta\sb M(t):=\sup\sb{s\in(0,t)}\abs{\zeta(s)}
=\zeta(t)$ for $0\le t<T\sb{\Omega}$,
and (\ref{bound-r}) takes the form
\begin{equation}\label{bound-r-2}
\abs{R\sb j(\eta,\zeta,\bm\upsilon)}
\le C\sb 6\zeta^2,
\qquad j=1,\,2,
\qquad
0\le t<T\sb{\Omega}.
\end{equation}
Using
(\ref{system-degenerate}) and (\ref{bound-r-2}),
and taking into account (\ref{zeta1-small})
and monotonicity of $\zeta(t)$ for $0\le t<T\sb{\Omega}$,
we compute:
\begin{equation}\label{dot-eta}
\dot\eta(t)=\zeta(t)+R\sb 1\ge \zeta(t)-C\sb 6\zeta^2(t)
=\zeta(t)(1-C\sb 6\zeta(t))>\zeta\sb 0/2,
\qquad
0\le t<T\sb{\Omega}.
\end{equation}
This allows to consider $\zeta$ as a function of $\eta$
(as long as $0\le t<T\sb{\Omega}$).
By (\ref{system-degenerate}), (\ref{bound-r-2}),
and (\ref{zeta1-small}),
\begin{equation}\label{dz-de}
\frac{d\zeta}{d\eta}
=\frac{\lambda(\eta)\zeta+R\sb 2}{\zeta+R\sb 1}
\le\frac{\lambda(\eta)\zeta+C\sb 6\zeta^2}{\zeta-C\sb 6\zeta^2}
=\frac{\lambda(\eta)+C\sb 6\zeta}{1-C\sb 6\zeta}
\le 2(\lambda(\eta)+C\sb 6\zeta),
\qquad
0\le t<T\sb{\Omega}.
\end{equation}
Thus,
$
\frac{d\zeta}{d\eta}-2C\sb 6\zeta<2\lambda(\eta)
$
for $0\le t<T\sb{\Omega}$.
Multiplying both sides of this relation by $e^{-2 C\sb 6\eta}$
and integrating, we get Gronwall's inequality:
\begin{equation}\label{non-gronwall}
\int\sb{\eta\sb 0}\sp{\eta}\frac{d}{d\eta'}
\left(e^{-2 C\sb 6\eta'}\zeta(\eta')\right)\,d\eta'
<2\int\sb{\eta\sb 0}\sp{\eta}e^{-2 C\sb 6\eta'}\lambda(\eta')\,d\eta'
\le 2 e^{-2 C\sb 6\eta\sb 0}\Lambda(\eta),
\end{equation}
\begin{equation}\label{asd}
\zeta
<e^{2 C\sb 6\eta}
\Bigl(
2 e^{-2 C\sb 6\eta\sb 0}\Lambda(\eta)
+
e^{-2 C\sb 6\eta\sb 0}\zeta\sb 0
\Bigr)
<3 e^{2 C\sb 6\eta}\Lambda(\eta),
\qquad
0\le t<T\sb{\Omega}.
\end{equation}
See Figure~\ref{fig-eta-zeta-mono}.
We used the inequality
$\zeta\sb 0<\Lambda(\eta\sb 0)\le\Lambda(\eta)$
that follows from (\ref{c-3})
and monotonicity of $\Lambda(\eta)$.

\begin{figure}[htbp]
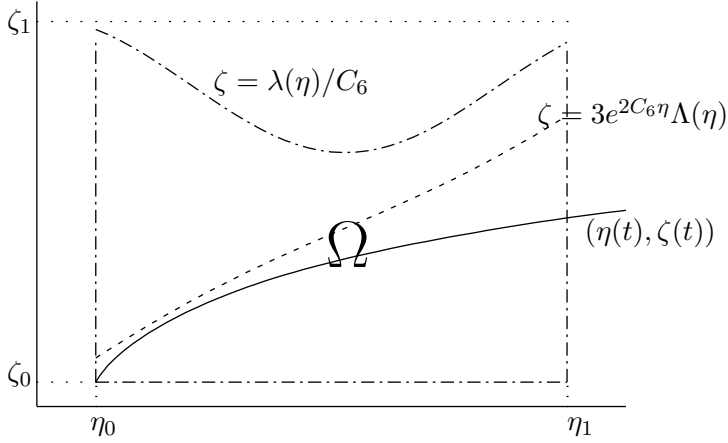

\begin{center}
\input kdvform-eta-zeta-mono.tex
\end{center}
\vskip -30pt
\caption{
The trajectory $(\eta(t),\zeta(t))$ (the solid line)
stays in the part of the region $\Omega$
below the dashed line $\zeta=3 e^{2C\sb 6\eta}\Lambda(\eta)$.
}
\label{fig-eta-zeta-mono}
\end{figure}

Now let us argue that $T\sb{\Omega}=T\sb 1$.
If $T\sb{\Omega}=\infty$, we are done,
therefore we only need to consider the case $T\sb{\Omega}<\infty$.
By (\ref{def-t22}),
the moment $T\sb{\Omega}$ is characterized by
\begin{equation}\label{t22-is-such}
{\it either}\quad
T\sb{\Omega}=T\sb 1
\quad{\it or}\quad
(\eta(T\sb{\Omega}),\zeta(T\sb{\Omega}))\in\p\Omega
\quad{\it or}\quad
\dot\eta(T\sb{\Omega})=0,
\end{equation}
or any combination of these three conditions.
By continuity,
the bound (\ref{asd}) is also valid at $T\sb{\Omega}$
(the last inequality in (\ref{asd}) remains strict);
therefore,
\begin{equation}\label{t2-3e}
\zeta(T\sb{\Omega})
< 3 e^{2 C\sb 6\eta(T\sb{\Omega})}\Lambda(\eta(T\sb{\Omega}))
<\lambda(\eta(T\sb{\Omega}))/C\sb 6.
\end{equation}
In the last inequality, we used Lemma~\ref{lemma-3e}.
The inequality (\ref{t2-3e}) also leads to
\begin{equation}\label{dot-zeta-t2}
\dot\zeta=\lambda(\eta)\zeta+R\sb 2
\ge\zeta(\lambda(\eta)-C\sb 6\zeta)>0,
\qquad
0\le t\le T\sb{\Omega}.
\end{equation}
Using (\ref{t2-3e}) and (\ref{dot-zeta-t2})
in (\ref{t22-is-such}),
we conclude that either $T\sb{\Omega}=T\sb 1$
or $\eta(T\sb{\Omega})=\eta\sb 1$ and
hence again $T\sb{\Omega}=T\sb 1$
(by (\ref{at-t3}), $\eta(t)<\eta\sb 1$ for $0\le t<T\sb 1$).
The bounds (\ref{monotone}) and (\ref{zeta-bounds})
for $0\le t<T\sb{\Omega}=T\sb 1$
follow from (\ref{dot-eta}) and (\ref{asd})
(note that $\dot\zeta\ge 0$ for $0\le t<T\sb{\Omega}=T\sb 1$
by (\ref{def-t22})).
\end{proof}

\begin{lemma}\label{lemma-o}
Assume that $\norm{\bm\rho\sb 0}\sb{H\sp 1}<\eta\sb 1$.
There exists $C\sb 7>0$
so that
\[
\norm{\bm\rho (t) }\sb{L\sb\mu\sp 2}
\le C\sb 7\Lambda(\eta),
\qquad
0\le t<T\sb 1.
\]
\end{lemma}

\begin{proof}
Using
the estimate (\ref{zeta-bounds})
from Lemma~\ref{lemma-zeta-le-Lambda-2}
and
the estimate (\ref{bound-v-mu})
from
Proposition~\ref{prop-5}
(where $\eta\sb M(t)=\eta(t)$ and $\zeta\sb M(t)=\zeta(t)$
due to (\ref{monotone}) and positivity of $\eta\sb 0$ and $\zeta\sb 0$),
we obtain:
\[
\norm{\bm\rho (t) }\sb{L\sb\mu\sp 2}
\le\abs{\zeta}\norm{{\bm e}\sb{3,c}}\sb{L\sb\mu\sp 2}
+\norm{\bm\upsilon}\sb{L\sb\mu\sp 2}
\le
\abs{\zeta}
\left(
\norm{{\bm e}\sb{3,c}}\sb{L\sb\mu\sp 2}
+
C\sb 4[\zeta+\mathscr{Y}(\eta)]
\right).
\]
Now the statement of the lemma follows from the bound (\ref{zeta-bounds}).
The value of $C\sb 7$ could be taken equal to $K(\eta\sb 1,\zeta\sb 1)$,
where
\begin{equation}\label{def-c7}
K(\eta\sb 1,\zeta\sb 1)
=3 e^{2 C\sb 6\eta\sb 1}
\left[
\sup\sb{c\in[c\sb\star,c\sb\star+\eta\sb 1]}
\norm{{\bm e}\sb{3,c}}\sb{L\sb\mu\sp 2}
+C\sb 4\zeta\sb 1
+
C\sb 4\left\{
2\eta\sb 1
+\abs{\mathscr{N}\sb{c\sb\star+\eta\sb 1}
-\mathscr{N}\sb{c\sb\star}}\sp{1/2}
\right\}
\right]
,
\end{equation}
where the term in the braces
dominates
$\mathscr{Y}(\eta)$ which was defined in (\ref{def-y}).
(When estimating $\mathscr{Y}(\eta)$
defined in (\ref{def-y}),
we used the bound $\norm{\bm\rho\sb 0}\sb{H\sp 1}<\eta\sb 1$.)
\end{proof}


\subsubsection*{Conclusion of the proof of Theorem~\ref{main-theorem}}
In Theorem~\ref{main-theorem},
let us take
\begin{equation}\label{def-epsilon}
\epsilon
=\min(
\fra{\kappa\eta\sb 1}{2},
\norm{\bm\phi\sb{c\sb\star}}\sb{H\sp 1}
)>0.
\end{equation}
Pick $\delta>0$ arbitrarily small.
To comply with the requirements of
Lemmas~\ref{lemma-init-data} and~\ref{lemma-o},
we may assume that $\delta$ is smaller than $\min(\eta\sb 1,\delta\sb 1)$.
Fix $\mu\in(0,\min(\mu\sb 0,\mu\sb 1))$,
with $\mu\sb 0$ from Assumption~\ref{assumption-spectrum}
and $\mu\sb 1$ as in Lemma~\ref{lemma-nd}.
Let $\eta\sb 0$ and $\zeta\sb 0$
satisfy all the inequalities
in Lemma~\ref{lemma-init-data};
then the conditions (\ref{init-well}) of Proposition~\ref{prop-well}
are satisfied.
Let
\[
{\bm u}\sb 0
=
\bm\phi\sb{c\sb\star+\eta\sb 0}
+\zeta\sb 0{\bm e}\sb{3,c\sb\star+\eta\sb 0},
\]
so that
${\bm u}\sb 0\in H\sp 2(\R)\cap L\sb{2\mu}\sp 2(\R)$
by (\ref{e-3-smooth})
and $\norm{{\bm u}\sb 0-\bm\phi\sb{c\sb\star}}\sb{H\sp 1}<\delta$
by (\ref{c-2}).
Proposition~\ref{prop-well} states that
there is $T\sb 1\in\R\sb{+}\cup\{+\infty\}$
and a function
${\bm u}(t)\in C([0,\infty),H\sp 2(\R)\cap L\sb{2\mu}\sp 2(\R))$,
${\bm u}(0)={\bm u}\sb 0$,
so that
for $0\le t<T\sb 1$
the function ${\bm u}(t)$ solves (\ref{gkdv})
and all the inequalities (\ref{at-t3})
are satisfied.

\begin{lemma}\label{lemma-finite}
In Proposition~\ref{prop-well},
one can only take $T\sb 1<\infty$.
\end{lemma}

\begin{proof}
If we had $T\sb 1=+\infty$,
then
$
\dot\eta
\ge\zeta\sb 0/2
$
for $t\in\R\sb{+}$
by Lemma~\ref{lemma-zeta-le-Lambda-2},
hence $\eta(t)$ would reach $\eta\sb 1$ in finite time,
contradicting the bound $\eta(t)<\eta\sb 1$
for $0\le t<T\sb 1$
from Proposition~\ref{prop-well}~({\it iii}).
\end{proof}

Since $T\sb 1<\infty$,
Proposition~\ref{prop-well}~({\it iv}) states that at least one of the inequalities
in (\ref{at-t3}) turns into equality at $t=T\sb 1$.
As follows from the bound (\ref{bound-v-mu})
and the inequality (\ref{cc-2}),
$\norm{\bm\upsilon(T\sb 1)}\sb{H\sb\mu\sp 1}<\delta\sb 1$.
Also, by (\ref{zeta-bounds})
(where the bound from above does not have to be strict at $T\sb 1$),
\begin{equation}\label{z-l-z1}
\zeta(T\sb 1)\le 3e^{2C\sb 6\eta\sb 1}\Lambda(\eta(T\sb 1))
\le 3e^{2C\sb 6\eta\sb 1}\Lambda(\eta\sb 1)
<\lambda(\eta)/C\sb 6
<\zeta\sb 1.
\end{equation}
We took into account the monotonicity of $\Lambda(\eta)$
and the inequalities (\ref{3e}) and (\ref{eta1-small}).
Therefore, either
$\norm{{\bm u}(T\sb 1)}\sb{H\sp 1}=2\norm{\bm\phi\sb{c\sb\star}}\sb{H\sp 1}$
or $\eta(T\sb 1)=\eta\sb 1$ (or both).
In the first case,
\begin{equation}\label{case-1}
\inf\sb{s\in\R}
\norm{{\bm u}(\cdot,T\sb 1)-\bm\phi\sb{c\sb\star}(\cdot-s)}\sb{H\sp 1}
\ge
\norm{{\bm u}(\cdot,T\sb 1)}\sb{H\sp 1}
-\norm{\bm\phi\sb{c\sb\star}}\sb{H\sp 1}
\ge\norm{\bm\phi\sb{c\sb\star}}\sb{H\sp 1}\ge\epsilon,
\end{equation}
hence the instability of $\bm\phi\sb{c\sb\star}$ follows.
We are left to consider the case $\eta(T\sb 1)=\eta\sb 1$.
According to (\ref{def-rho}),
\begin{eqnarray}
\inf\sb{s\in\R}\norm{{\bm u}(\cdot,t)-\bm\phi\sb{c\sb\star}(\cdot-s)}\sb{L\sp 2}
\ge
\inf\sb{s\in\R}\norm{{\bm u}(\cdot,t)
-\bm\phi\sb{c\sb\star}(\cdot-s)}\sb{L\sp 2(\R,\min(1,e\sp{\mu x})\,dx)}
\nonumber
\\
\ge
\inf\sb{s\in\R}\norm{\bm\phi\sb{c(t)}(\cdot)
-\bm\phi\sb{c\sb\star}(\cdot-s)}\sb{L\sp 2(\R,\min(1,e\sp{\mu x})\,dx)}
-\norm{\bm\rho(t)}\sb{L\sb\mu\sp 2}.
\label{ge2}
\end{eqnarray}
Applying Lemma~\ref{lemma-nd}
and
Lemma~\ref{lemma-o}
to the two terms in the right-hand side
of (\ref{ge2}),
we see that
\begin{equation}\label{case-2-0}
\inf\sb{s\in\R}\norm{{\bm u}(\cdot,t)
-\bm\phi\sb{c\sb\star}(\cdot-s)}\sb{L\sp 2}
\ge\kappa\eta-C\sb 7\Lambda(\eta),
\qquad
0\le t<T\sb 1,
\quad\kappa>0.
\end{equation}
Since 
$C\sb 7\Lambda(\eta\sb 1)\le\kappa\eta\sb 1/2$
by (\ref{c7-lambda}),
\begin{equation}\label{case-2}
\inf\sb{s\in\R}\norm{{\bm u}(\cdot,T\sb 1)
-\bm\phi\sb{c\sb\star}(\cdot-s)}\sb{L\sp 2}
\ge \kappa\eta\sb 1/2\ge\epsilon,
\end{equation}
and again the instability of $\bm\phi\sb{c\sb\star}$ follows.

This completes the proof
of Theorem~\ref{main-theorem}.

\section{Non-degenerate case: normal form}
\label{sect-normal-forms}

In this section, we prove that the critical soliton with the speed
$c\sb\star$ generally corresponds to the saddle-node bifurcation of
two branches of non-critical solitons. We assume for simplicity that
$c\sb\star$ is a non-degenerate critical point of
$\mathscr{N}\sb{c}$, in the sense that
\begin{equation}\label{npp-nz}
\mathscr{N}\sb{c\sb\star}'=0,
\qquad
\mathscr{N}\sb{c\sb\star}''\ne 0.
\end{equation}
We rewrite the two last equations from the system (\ref{xi-eta-zeta-system}):
\begin{equation}
\label{eta-zeta-system1}
\left[\begin{array}{cc}\dot{\eta}
\\
\dot{\zeta}\end{array}\right]
=\left[\begin{array}{ccc} 0 & 1
\\ 0 & \lambda\sb c\end{array}\right]
\left[\begin{array}{cc}\eta
\\
\zeta \end{array}\right]
+\left[\begin{array}{c}
R\sb 1(\eta,\zeta,\bm\upsilon)
\\ R\sb 2(\eta,\zeta,\bm\upsilon) \end{array}\right].
\end{equation}
As follows from (\ref{fourth}) and (\ref{def-lambda}),
\begin{equation}\label{Taylor-series-lambda}
\lambda\sb c =\lambda\sb{c\sb\star+\eta}
=\lambda\sb{c\sb\star}'\eta+\ol(\eta^2), \qquad
\lambda\sb{c\sb\star}' =-\frac{2
\mathscr{N}''\sb{c\sb\star}}{(I\sb{c\sb\star}')^2},
\end{equation}
where $ \lambda\sb{c\sb\star}'
\ne 0 $ by (\ref{npp-nz}). The system (\ref{eta-zeta-system1}) has
the nonlinear terms $R\sb j(\eta,\zeta,\bm\upsilon)$, $j=1,\,2$,
estimated in Lemma~\ref{lemma-bounds-on-r} for monotonically
increasing functions $\eta(t)$, $|\zeta(t)|$ on a local existence
interval $0 < t < T_1$. It follows from (\ref{def-r}) that
\[
R\sb 1(0,0,0)=R\sb 2(0,0,0)=0,
\]
so that the point $(\eta,\zeta)=(0,0)$ is a critical point of
(\ref{eta-zeta-system1}) when $\bm\upsilon=0$. This critical point
corresponds to the critical traveling wave
$\bm\phi\sb{c\sb\star}(x)$ itself. The following result establishes
a local equivalence between the system (\ref{eta-zeta-system1}) and
the truncated system $\ddot{\eta} = \lambda\sb{c\sb\star}' \eta \dot{\eta}$,
thus guaranteeing the instability of the critical point $(\eta,\zeta) = (0,0)$.

\begin{proposition}\label{prop-62}
Assume that the conditions (\ref{npp-nz}) are satisfied.
Consider the subset of trajectories
$(\eta(t),\zeta(t))$
of the system (\ref{eta-zeta-system1})
that lie inside the $\epsilon$-neighborhood
$\mathcal{D}\sb\epsilon\subset\R^2$ of the origin
and satisfy the condition that
both functions $\eta(t)$ and $|\zeta(t)|$ are monotonically
increasing.
For sufficiently small
$\epsilon>0$
this subset of the trajectories
is topologically equivalent
to a subset of the trajectories of the truncated normal form:
\begin{equation}
\label{normal-form3}\dot{x}=\frac{1}{2}\lambda\sb{c\sb\star}'
x^2+E\sb 1,
\end{equation}
where $E_1$ is constant.
\end{proposition}

\begin{proof}
Since $\zeta = \dot{\eta} - R_1(\eta,\zeta,\bm\upsilon)$, we can
rewrite the system (\ref{eta-zeta-system1}) in the equivalent form:
\begin{equation}
\label{untruncated-differential-form} \frac{d}{dt} \left( \dot{\eta}
- \frac{1}{2}\lambda\sb{c\sb\star}'  \eta^2 -
R_1(\eta,\zeta,\bm\upsilon) \right) = R(\eta,\zeta,\bm\upsilon),
\end{equation}
where
$$
R(\eta,\zeta,\bm\upsilon)
\equiv R_2(\eta,\zeta,\bm\upsilon) - \lambda_c
R_1(\eta,\zeta,\bm\upsilon)) + (\lambda_c - \lambda\sb{c\sb\star}'
\eta )\zeta.
$$
It follows from Lemma~\ref{lemma-bounds-on-r} and
(\ref{Taylor-series-lambda}) that there exists a constant $C > 0$
such that $|R| \leq C (\zeta^2 + \eta^2 |\zeta|)$. The integral form
of (\ref{untruncated-differential-form}) is
\begin{equation}
\label{untruncated-normal-form} \dot{\eta} -
\frac{1}{2}\lambda\sb{c\sb\star}'  \eta^2 - E_1  =
\tilde{R}(t),
\end{equation}
where
$$
\tilde{R}(t)\equiv R_1(\eta(t),\zeta(t),\bm\upsilon(t))
+\int_0^t R(\eta(t'),\zeta(t'),\bm\upsilon(t'))\,dt'
$$
and $E_1$ is the constant of integration.
Using Lemma~\ref{lemma-bounds-on-r},
the bound $|\zeta| \leq \dot{\eta} +C_6 \zeta^2$,
and integration by parts, we obtain that
\[
\int_0^t \zeta^2\,dt'
\le \eta\abs{\zeta}+ C_6 \int_0^t\abs{\zeta}^3\,dt'
\le\eta\abs{\zeta}+ C_6\eta\abs{\zeta}^2
+ C_6^2 \int_0^t\abs{\zeta}^4\,dt'
\le
\dots
\le
\frac{\eta\abs{\zeta}}{1 - C_6\abs{\zeta}}
\]
and
\[
\int_0^t \eta^2 |\zeta|\,dt'
\le\frac{\eta^3}{3}+ C_6 \int_0^t\eta^2\zeta^2\,dt'
\le\dots\le\frac{\eta^3}{3(1 - C_6 |\zeta|)}.
\]
Thus, if $\abs{\zeta}$ is sufficiently small, there exists a constant
$\tilde{C}>0$
such that $|\tilde{R}|\le\tilde{C}(\zeta^2+|\zeta|\eta+\eta^3)$.
The topological equivalence of equation (\ref{untruncated-normal-form})
with the above estimate on
$|\tilde{R}|$ in the disk $(\eta,\zeta) \in
\mathcal{D}\sb{\epsilon}$ to the truncated normal form
(\ref{normal-form3}) with sufficiently small $E_1$ is proved in
\cite[Lemma 3.1]{MR1711790}.
By definition, two systems are said to be topologically equivalent if
there exists a homeomorphism between solutions of these systems.
We note that this equivalence holds for a family of trajectories which
corresponds to monotonically increasing functions $\eta(t)$,
$|\zeta(t)|$ in a subset of the small disk near $(\eta,\zeta) = (0,0)$.
\end{proof}

\begin{corollary}
The critical point $(0,0)$ of system (\ref{eta-zeta-system1})
is unstable, in the sense that
there exists $\epsilon>0$ such that for any $\delta>0$
there are $(\eta(0),\zeta(0))\in\mathcal{D}\sb{\delta}$
and $t\sb\ast=t\sb\ast(\delta,\epsilon)<\infty$
such that $(\eta(t\sb\ast),\zeta(t\sb\ast))\notin\mathcal{D}\sb{\epsilon}$.
\end{corollary}

\begin{proof}
The normal form equation (\ref{normal-form3}) shows that the
critical point $x = 0$ is semi-stable at $E_1 = 0$, such that the
trajectory with any $x(0) \neq 0$ of the same sign as
$\lambda\sb{c\sb\star}'$ escapes the local neighborhood of the point
$x = 0$ in a local time $t \in [0,T]$. By Proposition~\ref{prop-62},
local dynamics of (\ref{normal-form3}) for $x(t)$ is equivalent to
local dynamics of (\ref{eta-zeta-system1}) for $(\eta,\zeta)$.
\end{proof}

\begin{remark}
The truncated normal form (\ref{normal-form3}) is rewritten for
$c=c\sb\star+x$:
\begin{equation}
\label{normal-form4}\dot{c}=\frac{1}{2}\lambda\sb{c\sb\star}'(c-c\sb\star)^2
+E\sb 1.
\end{equation}
The normal form (\ref{normal-form4}) corresponds to the standard
saddle-node bifurcation.
It was derived and studied in \cite{MR1416294}
by using the asymptotic multi-scale expansion method.
When $E=0$, the critical point $c=c\sb\star$ is a degenerate saddle point,
which is nonlinearly unstable.
Assume for definiteness that
$\lambda\sb{c\sb\star}'>0$
(which implies that $\mathscr{N}''\sb{c\sb\star}<0$).
Then there are no fixed points for $E\sb 1>0$
and two fixed points for $E\sb 1<0$
in the normal form equation (\ref{normal-form4}).
Therefore,
there exist initial perturbations (with $E\sb 1>0$ and any $c\sb 0$
or with $E\sb 1=0$ and $c\sb 0>c\sb\star$),
which are arbitrarily close to
the traveling wave with $c=c\sb\star$, but the norm $\abs{c-c\sb\star}$
exceeds
some a priori fixed value
at $t=t\sb\ast>0$.
Two fixed points exist
for $E\sb 1<0$:
\begin{equation}
\label{fixed-point} c=c\sb{E}\sp{\pm}=c\sb\star
\pm \sqrt{\frac{E\sb 1}{\mathscr{N}''\sb{c\sb\star}}}\abs{I\sb{c\sb\star}'},
\end{equation}
so that $c=c\sb{E}\sp{+}$ is an unstable saddle point
and $c=c\sb{E}\sp{-}$ is a stable node.
The two fixed points correspond to two branches of
traveling waves with $\mathscr{N}\sb{c}<\mathscr{N}\sb{\rm max}$,
where $\mathscr{N}\sb{\rm max}=\mathscr{N}(\bm\phi\sb{c\sb\star})$.
The left branch with $c\sb{E}\sp{-} < c\sb\star$
corresponds to $\mathscr{N}'\sb{c\sb{E}\sp{-}}>0$ and the right branch with
$c\sb{E}\sp{+}>c\sb\star$
corresponds to $\mathscr{N}'\sb{c\sb{E}\sp{+}}<0$.
According to the stability theory for traveling waves \cite{MR1177566},
the left branch is orbitally stable,
while the right branch is linearly unstable.
\end{remark}

\appendix

\section{Appendix: Existence of solitary waves}\label{sect-existence}

Let us discuss the existence of standing waves.
We assume that $f$ is smooth.
Let $F$ denote the primitive of $f$
such that $F(0)=0$.
Thus, by (\ref{f-p-zero}),
\begin{equation}\label{f-is-zero}
F(0)=F'(0)=F''(0)=0.
\end{equation}
The wave profile
$\bm\phi\sb{c}$ is to satisfy the equation
\[
u''-c u
=f(u),
\qquad c>0.
\]
Multiplying this by $u'$
and integrating,
and taking into account
that we need
$\lim\sb{\abs{x}\to\infty}u(x)=0$,
we get
\begin{equation}\label{phi-prime}
\frac{d u(x)}{dx}
=\pm\sqrt{c u^2+2 F(u)}.
\end{equation}
There will be a strictly positive continuous solution
exponentially decaying at infinity
if there exists $\bm\xi\sb{c}>0$ such that
$c\frac{u^2}{2}+F(u)>0$
for $0<u<\bm\xi\sb{c}$,
and also
\[
c\frac{\bm\xi\sb{c}^2}{2}+F(\bm\xi\sb{c})=0,
\qquad
c\bm\xi\sb{c}+f(\bm\xi\sb{c})<0.
\]
The last two conditions imply that the map $c\mapsto\bm\xi\sb{c}$
is invertible and smooth (as $F$ is).
One immediately sees that $\bm\phi\sb c\in C\sp\infty(\R)$
and, due to the exponential decay at infinity,
$\bm\phi\sb c\in H\sp\infty(\R)$.
For each $c$, the solution $\bm\phi\sb c$ is unique
(up to translations of the origin),
and (after a suitable translation of the origin)
satisfies the following properties:
it is strictly positive, symmetric,
and is monotonically
decreasing (strictly) away from the origin.
This result
follows from the implicit representation
\begin{equation}\label{x-vs-phi}
x=\pm\int\sb{\bm\phi\sb c}\sp{\bm\xi\sb{c}}
\frac{du}{\sqrt{c u^2+2 F(u)}}.
\end{equation}
See \cite[Section 6]{MR695535}
for the exhaustive treatment of this subject.


\begin{lemma}\label{lemma-exp}
There exist positive constants
$C\sb 1$, $C\sb 2$, $C\sb 1'$, and $C\sb 2'$ such that
\begin{equation}\label{bounds-on-phi}
C\sb 1 e^{-\sqrt{c}\abs{x}}
\le \abs{\bm\phi\sb{c}(x)}
\le C\sb 2 e^{-\sqrt{c}\abs{x}},
\qquad x\in\R,
\end{equation}
\begin{equation}\label{bounds-on-e1}
C\sb 1' e^{-\sqrt{c}\abs{x}}
\le \abs{\p\sb x\bm\phi\sb{c}(x)}
\le C\sb 2' e^{-\sqrt{c}\abs{x}},
\qquad \abs{x}\ge 1.
\end{equation}
\end{lemma}

\begin{proof}
Since $\lim\sb{\abs{x}\to\infty}\bm\phi\sb{c}(x)=0$,
there exists $x\sb 1>0$ so that
$\abs{F(\bm\phi\sb{c}(x))}/\bm\phi\sb{c}^2(x)<c/4$
for $\abs{x}\ge x\sb 1$.
Then, for $x>x\sb 1$,
we get from (\ref{x-vs-phi}):
\[
x-x\sb 1
=
\int\sb{\bm\phi\sb{c}(x)}\sp{\bm\phi\sb{c}(x\sb 1)}
\frac{d u}{\sqrt{c u^2+2 F(u)}}.
\]
It follows that
\begin{equation}\label{the-above}
\int\sb{\bm\phi\sb{c}(x)}\sp{\bm\phi\sb{c}(x\sb 1)}
\frac{d u}{c\sp{1/2}u}
-
\int\sb{\bm\phi\sb{c}(x)}\sp{\bm\phi\sb{c}(x\sb 1)}
\frac{\abs{F(u)}}{c\sp{3/2}u^3}\,du
\le
x-x\sb 1
\le
\int\sb{\bm\phi\sb{c}(x)}\sp{\bm\phi\sb{c}(x\sb 1)}
\frac{d u}{c\sp{1/2}u}
+
\int\sb{\bm\phi\sb{c}(x)}\sp{\bm\phi\sb{c}(x\sb 1)}
\frac{\abs{F(u)}}{c\sp{3/2}u^3}\,du.
\end{equation}
By (\ref{f-is-zero}),
$\abs{F(u)}/u^3$ is bounded for $u$ small,
and we conclude from (\ref{the-above}) that
\begin{equation}\label{bound-phi}
\ln\bm\phi\sb{c}(x)
-C\sb 3
\le
c\sp{1/2}(x-x\sb 1)
\le
\ln\bm\phi\sb{c}(x)+C\sb 3,
\end{equation}
where
$
C\sb 3
=
c\sp{-1}
\int\sb 0\sp{\bm\phi\sb{c}(x\sb 1)}
\abs{F(u)}u^{-3}\,du.
$
Inequalities (\ref{bound-phi})
immediately prove (\ref{bounds-on-phi}).
Bounds (\ref{bounds-on-e1})
immediately follow from (\ref{phi-prime}).
\end{proof}

We also need the following result
that gives the rate of decay of 
${\bm e}\sb{2,c}=\p\sb c\bm\phi\sb c$
and ${\bm e}\sb{3,c\sb\star}$ at infinity.

\begin{lemma}\label{lemma-decay}
Let $R\in C\sp\infty(\R)$
satisfy the bound
$\abs{R(x)}\le C\sb 1 e^{-\sqrt{c}\abs{x}}$
for $x \ge 0$,
for some $c>0$, $C\sb 1>0$.
Let $u\in C\sp\infty(\R)$ satisfy
\begin{equation}\label{u-cu-r}
u''-c u=R,
\qquad
\lim\sb{x\to +\infty}u(x)=0.
\end{equation}
Then there exists $C\sb 2>0$ (that depends on $c$, $C\sb 1$, and $u$)
such that
\begin{equation}\label{ubbb}
\abs{u(x)}\le C\sb 2(1+\abs{x})e^{-\sqrt{c}\abs{x}},
\qquad x\ge 0.
\end{equation}
\end{lemma}

\begin{remark}
$C\sb 2$ depends not only on $c$ and $C\sb 1$ but also on $u$
because the solution to (\ref{u-cu-r})
is defined up to $\const e^{-\sqrt{c}\,x}$.
\end{remark}

\begin{proof}
First, we notice that if
$P\in C\sp\infty(\R)$, $P(x)\ge 0$ for $x\ge 0$,
and if $v\in C\sp\infty(\R)$ solves
\begin{equation}\label{m-p}
v''-c v=P(x),\qquad
v(0)=0,\qquad \lim\sb{x\to +\infty}v(x)=0,
\end{equation}
then $v(x)\le 0$ for $x\ge 0$.
(The existence of a point $x\sb 0>0$
where $u$ assumes a positive maximum
contradicts the equation in (\ref{m-p}).)

Now we consider the functions $u\sb{-}$ and $u\sb{+}$
that satisfy
\begin{equation}
u\sb\pm''(x)-c u\sb\pm=\pm C\sb 1 e^{-\sqrt{c}\abs{x}},
\qquad
u\sb\pm(0)=u(0),
\qquad\lim\sb{x\to +\infty}u\sb\pm(x)=0.
\end{equation}
Both $u\sb\pm$ can be written explicitly;
they satisfy (\ref{ubbb}).
Since $v=u-u\sb{-}$ and $v=u\sb{+}-u$
satisfy (\ref{m-p}) with 
$P(x)=C\sb 1 e^{-\sqrt{c}\abs{x}}+R(x)$
and $P(x)=C\sb 1 e^{-\sqrt{c}\abs{x}}-R(x)$,
respectively,
we conclude that
$u\sb{+}(x)\le u(x)\le u\sb{-}(x)$ for $x\ge 0$,
hence $u$ also satisfies (\ref{ubbb}).
\end{proof}

\section{Appendix: Fredholm alternative for $\mathcal{H}\sb{c}$}
\label{sect-sd}

\begin{lemma}[Fredholm Alternative]\label{lemma-Fredholm}
Let $R(x)\in\mathscr{S}\sb{+,m}(\R)$,
$m\ge 0$
(see Definition~\ref{def-Schwarz}).
If
\begin{equation}
\int\sb{\R}{\bm e}\sb{1,c}(x)R(x)\,dx=0,
\end{equation}
then the equation
\begin{equation}\label{hu-r}
\mathcal{H}\sb{c}u=R
\end{equation}
has a solution $u\in\mathscr{S}\sb{+,m}(\R)$.
(This solution is unique if we impose the constraint
$\langle {\bm e}\sb{1,c},u\rangle=0$.)
Otherwise, any solution
$u(x)$ to (\ref{hu-r})
such that $\lim\sb{x\to+\infty}u(x)=0$
grows exponentially at $-\infty$:
\[
\lim\sb{x\to-\infty}e^{-\sqrt{c}\abs{x}}u(x)
\ne 0.
\]
\end{lemma}

\begin{proof}
Let us pick an even function $R\sb{+}\in H\sp\infty(\R)$
so that $R\sb{+}(x)=R(x)$ for $x\ge 1$.
Since $R\sb{+}$ is even
and therefore orthogonal to the kernel of the operator
$\mathcal{H}\sb{c}$, there is a solution
$u\sb{+}\in H\sp\infty(\R)$ to the equation
\begin{equation}
\mathcal{H}\sb{c}u\sb{+}=R\sb{+}.
\end{equation}
Denote by $u$
the solution to the ordinary differential equation
\begin{equation}\label{hf-phi}
\mathcal{H}\sb{c}u
\equiv
-u''+(f'(\bm\phi\sb{c})+c)u=R
\end{equation}
such that
$u\at{x=1}=u\sb{+}\at{x=1}$,
$u'\at{x=1}=u\sb{+}'\at{x=1}$.
Then $u\in C\sp\infty(\R)$
coincides with $u\sb{+}$
for $x\ge 1$
and thus satisfies
\begin{equation}\label{lim-f-infty}
\lim\sb{x\to+\infty}u(x)=0.
\end{equation}
We take the pairing of (\ref{hf-phi}) with
${\bm e}\sb{1,c}$:
\begin{equation}\label{e1-hf-phi}
\int\sb{x}\sp{\infty}
{\bm e}\sb{1,c}(y)\mathcal{H}\sb{c}u(y)\,dy
=\int\sb{x}\sp{\infty}{\bm e}\sb{1,c}(y)R(y)\,dy
\equiv r(x),
\qquad
x\in\R.
\end{equation}
Since
\[
{\bm e}\sb{1,c}\mathcal{H}\sb{c}u
=
u
\mathcal{H}\sb{c}{\bm e}\sb{1,c}
-{\bm e}\sb{1,c}\p\sb{x}^2u
+u\p\sb{x}^2{\bm e}\sb{1,c}
=
-\p\sb{x}({\bm e}\sb{1,c}u')
+\p\sb{x}(u\p\sb{x}{\bm e}\sb{1,c}),
\]
where we took into account that
$\mathcal{H}\sb{c}{\bm e}\sb{1,c}=0$,
we obtain from (\ref{e1-hf-phi}) the relation
\begin{equation}\label{f-ode}
{\bm e}\sb{1,c}(x)u'(x)
-u(x)\p\sb{x}{\bm e}\sb{1,c}(x)
=r(x).
\end{equation}
The boundary term at $x=+\infty$
does not contribute into (\ref{f-ode})
due to the limit (\ref{lim-f-infty}).
We will use this relation to find the behavior of $u(x)$
as $x\to-\infty$.
For $x\le-1$,
we divide
the relation (\ref{f-ode})
by ${\bm e}\sb{1,c}^2$
(we can do this since ${\bm e}\sb{1,c}(x)=-\p\sb{x}\bm\phi\sb{c}(x)\ne 0$
for $x\ne 0$),
getting
\begin{equation}\label{f-ode-2}
\p\sb{x}
\left(
\frac{u(x)}{{\bm e}\sb{1,c}(x)}
\right)
=\frac{r(x)}{{\bm e}\sb{1,c}^2(x)}.
\end{equation}
Therefore, for $x\le -1$,
\begin{equation}\label{f-is}
u(x)
-
{\bm e}\sb{1,c}(x)\frac{u(-1)}
{{\bm e}\sb{1,c}(-1)}
={\bm e}\sb{1,c}(x)
\int\sb{-1}\sp{x}
\frac{r(y)\,dy}{{\bm e}\sb{1,c}^2(y)}
={\bm e}\sb{1,c}(x)
\int\sb{-1}\sp{x}
\frac{(r(y)-r\sb{-})+r\sb{-}}{{\bm e}\sb{1,c}^2(y)}\,dy,
\end{equation}
where $r\sb{-}=\lim\sb{x\to-\infty}r(x)$.

Since $R\in\mathscr{S}\sb{+,m}(\R)$,
$\abs{R(x)}\le C(1+\abs{x})\sp m$,
$m\in\Z$,
$m\ge 0$.
Using Lemma~\ref{lemma-exp},
we see that
\begin{equation}\label{q-q-bound}
\Abs{r(x)-r\sb{-}}
=\Abs{\int\sb{-\infty}\sp{x}R(y){\bm e}\sb{1,c}(y)\,dy}
\le \const e^{-\sqrt{c}\abs{x}}(1+\abs{x})\sp m,
\qquad
x\le -1.
\end{equation}
At the same time,
Lemma~\ref{lemma-exp}
also shows that
\begin{equation}
\int\sb{-1}\sp{x}\frac{dy}{{\bm e}\sb{1,c}^2(y)}
\ge\const e^{2\sqrt{c}\abs{x}},
\qquad
x\le -1.
\end{equation}
Therefore, if $r\sb{-}\ne 0$,
the right-hand side of (\ref{f-is})
grows exponentially
as $x\to-\infty$.
The same is true for $u(x)$,
since
the second term in the left-hand side of (\ref{f-is})
decays exponentially when $\abs{x}\to\infty$
by Lemma~\ref{lemma-exp}.
If instead $r\sb{-}=0$,
Lemma~\ref{lemma-exp}
and
the bound (\ref{q-q-bound})
show that
the right-hand side of (\ref{f-is})
is bounded by $\const (1+\abs{x})\sp m$,
proving similar bound for $u(x)$.
Using (\ref{hf-phi}) to get the bounds
on the derivatives
$u\sp{(N)}$, we conclude that
$u\in\mathscr{S}\sb{+,m}(\R)$.
\end{proof}

\section{Appendix: non-degeneracy of
$
\inf\limits\sb{s\in\R}\norm{\bm\phi\sb{c}(\cdot)
-\bm\phi\sb{c\sb\star}(\cdot-s)}
$
at $c\sb\star$
}
\label{sect-nd}

\begin{lemma}\label{lemma-nd}
If $\eta\sb 1>0$ is sufficiently small,
there exist $\mu\sb 1>0$ and $\kappa>0$
so that
\[
\inf\sb{s\in\R}\norm{\bm\phi\sb{c}(\cdot)
-\bm\phi\sb{c\sb\star}(\cdot-s)}\sb{L\sp 2(\R,\min(1,e\sp{\mu x})\,dx)}
\ge\kappa\abs{c-c\sb\star},
\qquad
c\in[c\sb\star,c\sb\star+\eta\sb 1],
\qquad
\mu\in[0,\mu\sb 1].
\]
\end{lemma}

\begin{proof}
Consider the function
\begin{equation}
g\sb\mu(c,s)
=
\norm{\bm\phi\sb{c}(\cdot)-\bm\phi\sb{c\sb\star}(\cdot-s)}
\sb{L\sp 2(\R,\min(1,e\sp{\mu x})\,dx)}^2.
\end{equation}
It is a smooth non-negative
function of $c$ and $s$, for $c\in[c\sb\star,c\sb\star+\eta\sb 1]$
and $s\in\R$.
It also depends smoothly on the parameter $\mu\ge 0$.
Zero is its absolute minimum, achieved at the point
$(c,s)=(c\sb\star,0)$.
We also note that the point $(c\sb\star,0)$
is non-degenerate when $\mu=0$:
\[
\p\sb{c}^2 g\sb 0(c,s)\at{(c\sb\star,0)}
=2\norm{\p\sb{c}\bm\phi\sb{c}\at{c=c\sb\star}}\sb{L\sp 2}^2
>0,
\qquad
\p\sb s^2 g\sb 0(c,s)\at{(c\sb\star,0)}
=2\norm{\p\sb{x}\bm\phi\sb{c\sb\star}}\sb{L\sp 2}^2
>0,
\]
\[
\p\sb{c}\p\sb s g\sb 0(c,s)\at{(c\sb\star,0)}
=-2(\p\sb{c}\bm\phi\sb{c}\at{c=c\sb\star},\p\sb{x}\bm\phi\sb{c\sb\star})
=0.
\]
By continuity, the quadratic form $g\sb\mu''\at{(c\sb\star,0)}$
is non-degenerate for $0\le \mu\le\mu\sb 1$,
with some $\mu\sb 1>0$.
Therefore, there exists $\kappa>0$
and an open neighborhood
$\Omega\subset\R^2$
of the point $(c\sb\star,0)$
such that
\begin{equation}\label{g-mu-non-degenerate}
g\sb\mu(c,s)\ge\kappa^2((c-c\sb\star)^2+s^2),
\qquad
(c,s)\in\Omega,
\qquad
0\le\mu\le\mu\sb 1.
\end{equation}
Moreover, we claim that
\begin{equation}\label{ge-epsilon}
\varGamma
\equiv
{\inf\sb{\mu\in(0,\mu\sb 1)}}
\;\;
{\inf\sb{(c,s)\in[c\sb\star,c\sb\star+\eta\sb 1]\times\R)\backslash\Omega}}
g\sb\mu(c,s)
>0.
\end{equation}
To prove (\ref{ge-epsilon}),
we only need to note that $(c\sb\star,0)$
is the only point where $g\sb\mu(c,s)$ takes the zero value
and that
$\lim\sb{\abs{s}\to\infty}g\sb\mu(c,s)
\ge
\inf\sb{c\in[c\sb\star,c\sb\star+\eta\sb 1]}
\norm{\bm\phi\sb{c}}\sb{L\sp 2(\R,\min(1,e\sp{\mu\sb 1 x})\,dx)}^2
>0$.

Now, we assume that $\eta\sb 1>0$
is small enough so that $\kappa^2\eta\sb 1^2<\varGamma$.
Then, by (\ref{g-mu-non-degenerate})
(valid for $(c,s)\in\Omega$)
and (\ref{ge-epsilon})
(valid for $(c,s)\in([c\sb\star,c\sb\star+\eta\sb 1]\times\R)\backslash\Omega$),
we conclude that
\begin{equation}
\inf\sb{s\in\R}g\sb\mu(c,s)\ge\kappa^2(c-c\sb\star)^2,
\qquad
c\in[c\sb\star,c\sb\star+\eta\sb 1],
\quad
\mu\in[0,\mu\sb 1].
\end{equation}
This proves the Lemma.
\end{proof}

\bibliographystyle{amsalpha}
\bibliography{all,books,strauss,kdvform,kdvform-local}

\end{document}